\documentclass{amsart}%
\usepackage{amssymb}
\usepackage{amsfonts}
\usepackage{amsmath}
\usepackage{comment}
\usepackage{mathtools}
\usepackage{mathrsfs}
\usepackage{graphicx}%
\usepackage{cite}

\newtheorem{theorem}{Theorem}
\theoremstyle{plain}

\newtheorem{corollary}{Corollary}

\newtheorem{definition}{Definition}
\newtheorem{example}{Example}

\newtheorem{lemma}{Lemma}

\newtheorem{remark}{Remark}

\numberwithin{equation}{section}
\begin{document}
\title[ ]{Neumann Data and Second Variation Formula of Renormalized Area for Conformally Compact Static Spaces }
\author{Zhixin Wang}
\address{Department of Mathematics, Shanghai Jiao Tong University, Shanghai, 201100}
\email{jhin@sjtu.edu.cn}

\begin{abstract}
In this paper, we derive the first and second variation formulas for the renormalized area for static Einstein spaces along a specific direction, demonstrating that the negativity of the Neumann data implies instability. Consequently, we obtain a rigidity result for the case when the conformal boundary is a warped product torus, which strengthens the result presented in \cite{GSW}.
\end{abstract}

\maketitle

\section{Introduction}
General relativity, proposed by Einstein in 1915, revolutionized our understanding of gravity and provided a new theoretical framework for exploring the large-scale structure of the universe. Within this framework, static Einstein manifolds play a fundamental role. Specifically, given a Riemannian manifold $(M,g)$, if there exists a smooth function $V$ that satisfies the equation
\begin{equation}
\nabla^2 V - \Delta V g - V \text{Ric} = 0 \label{static_Einstein_equ}
\end{equation} 
then $(M,g,V)$ is called a (vacuum) static Einstein manifold, and the equation is referred to as the static Einstein equation, with $V$ being called the potential.

The study of static Einstein manifolds has greatly enhanced our understanding of curvature properties. Firstly, by \cite{besse2007einstein} $(M,g,V)$ satisfying (\ref{static_Einstein_equ}) can be lifted to get manifolds of constant Ricci curvature through
\begin{equation}
(M\times \mathbb{R},\hat{g}_{\pm} = \pm V^2 dt^2 + g) \label{lifting}
\end{equation}

Secondly, the left-hand side of equation \eqref{static_Einstein_equ}, as the conjugate operator of the scalar curvature linearized operator, plays a crucial role in scalar curvature studies, especially in prescribed curvature problems \cite{Bo,Co,FM}.

It's well known that static Einstein spaces has constatn scalar curvature. Based on the sign of the scalar curvature, static Einstein manifolds can be classified into three types:

1)$R>0$, represented by the sphere; 

2)$R=0$, represented by Euclidean space and the Schwarzschild solution;

3)$R<0$, represented by hyperbolic space and the Horowitz-Myers soliton.
Among these, the third type generally has a conformally compact property, and it is the main focus of this study.
\begin{definition}\label{def_APE}
    A non-compact complete Riemannian manifold $(M^n,g)$ is called conformally compactifiable if $M$ is the interior of a compact manifold $\bar{M}$ with non-empty boundary $\Sigma$ and there exists a boundary defining function $x$ satisfying $x=0$ on $\Sigma$, $|dx|\neq 0$ on $\Sigma$, $x>0$ in $\mathring{M}$ so that $x^2g$ extends to a smooth metric on $\bar{M}$. $(\Sigma, x^2g\big|_{\Sigma})$ is called the conformal boundary. Furthermore,
    \begin{enumerate}
        \item $(M,g)$ is called to be Poincar\'e-Einstein (PE) manifold if $Ric=-(n-1)g$;
        \item $(M,g)$ is called to be asymptotically-Poincar\'e-Einstein (APE) if $Ric=-(n-1)g+O(x^n)$;
        \item $(M,g)$ is called to satisfy \textbf{Condition C} if there exists a $V$ so that $\frac{1}{V}$ is a bdf and $(M,g,V)$ satisfies (\ref{static_Einstein_equ}). 
    \end{enumerate}
\end{definition}

In APE case, one can choose different boundary defining functions (bdf's), each of which induces a conformal change in the boundary metric \((\Sigma, g)\). Therefore, it is more natural to consider the conformal class of the boundary metric rather than a specific representative. 

However, if \((M, g, V)\) satisfies \textbf{Condition C}, then \(\frac{1}{V}\) provides a natural and canonical choice of bdf. As a result, we obtain a distinguished metric within the conformal class, making this choice intrinsic to the geometric structure.\\

The classification of static Einstein manifolds has always been an important topic and has been explored in different contexts. When $R>0$, if $(M,g)$ is locally conformally flat, a general classification result was obtained in \cite{Ko,La}. While for $R=0$, especially asymptotically flat case, the rigidity is closely related to No-Hair theorem\cite{israel1968event}.

In the conformal compact case($R<0$), fruitful results were obtained for different conformal boundaries. The study of the rigidity of conformally compact static Einstein manifolds dates back to \cite{boucher1984uniqueness}, where the uniqueness of static Einstein manifolds with the standard sphere as the conformal boundary in three dimensions was proven, under the assumption of the positive mass conjecture for asymptotically hyperbolic manifolds. Later, X.Wang further improved this result, proving rigidity in the cases where $(M,g)$ has a spin structure or $n \leq 7$ \cite{wang2005uniqueness}. In the case where $M$ has a spin structure, the proof relies on an integral identity, which allows the detection of the Wang mass of asymptotically hyperbolic manifolds. Using the positive mass theorem for asymptotically hyperbolic space \cite{Wa}, part of the rigidity problem was solved. For $n \leq 7$, inspired by \cite{qing2003rigidity}, by choosing $\frac{1}{V+1}$ instead of $\frac{1}{V}$ as the boundary defining function, it was shown that $(M, \frac{1}{(V+1)^2}g)$ has non-negative scalar curvature, with the boundary being the standard sphere and the mean curvature being $n-1$. By combining the positive mass theorem for asymptotically flat manifolds in the case of $n \leq 7$, rigidity was proven.

When the conformal boundary is flat torus, we have the Horowitz-Myers soliton given by

         \begin{equation}
         \begin{aligned}
         M&=\mathbb{R}^2\times \mathbb{T}^{n-2}\\
             g&=\frac{1}{r^2(1-\frac{1}{r^n})}dr^2+r^2(1-\frac{1}{r^n})d\theta^2+r^2\sum_{i=1}^{n-2}ds_i^2\\
             V&=r\label{Horowitz-Myers}
             \end{aligned}
         \end{equation}

where $r\in [ 1,\infty),\theta\in [0,\frac{4\pi}{n}), s\in [0,a_i]$ with identification $\theta \sim \theta+\frac{4\pi}{n}$ and $s_i\sim s_i+a_i$. $a_i$'s are arbitrary positive numbers. When the conformal boundary satisfies certain convex condition, minimizing geodesic line  can be constructed in the universal covering of the compactification $(M^n, \frac{1}{V^2}g)$, and this geodesic line can be lifted to a non-chronal null line in the Lorentz metric obtained via (\ref{lifting}). By combining the maximum principle for Lorentz metrics\cite{Ga1}, a rigidity result was obtained in \cite{GSW}. Later in \cite{wang2024riccati}, this result was reproved using Busemann function without involving knowledge from general relativity. However, both methods require constructing geodesic lines near the conformal boundary, which imposes certain convexity conditions near the boundary of $\frac{1}{V^2} g$. Whether this condition can be removed remains open.

The study of Poincaré-Einstein manifolds runs parallel to the study of conformally compact static Einstein manifolds, with rigidity results in this area also attracting significant attention. Especially in the asymptotically hyperbolic case, a series of works \cite{qing2003rigidity} \cite{shi2005rigidity} address this, and the issue was ultimately resolved in \cite{DJ}. In \cite{DJ}, the authors used the Bishop volume comparison theorem and the Yamabe invariant of the sphere.\\

This paper is motivited by the following observation. A key property of the static Einstein triple $(M^n, g, V)$ is that it can be lifted to an $(n+1)$-dimensional Lorentz spacetime $(M\times \mathbb{R},\hat{g}_{-})$ of constant Ricci curvature as in (\ref{lifting}). In particular, $Ric_{\hat{g}_{-}}(\nu,\nu)=0$ for null vector $\nu$. Projecting this structure back to $(M,g,V)$, we find that the static Einstein triple inherits certain features of Ricci-flat geometry along a specific geometric flow. More precisely, consider the flow
\begin{equation}
\begin{aligned}
    \Phi: F\times I\rightarrow M\\
    \frac{\partial \Phi}{\partial t}=V \nu\label{flow}
\end{aligned}
\end{equation}
where $I$ is some open interval, $V$ is the static potential and $\nu$ is the unit normal vector. Let $F_t=\{\Phi(\cdot,t)\}$. Along this flow, we have the inequality
\begin{equation}
    \frac{1}{V}\frac{\partial }{\partial t}\frac{H}{V}\leq -\frac{1}{n-1}(\frac{H}{V})^2\label{Riccati_equ}
\end{equation}
which parallels the inequality
\begin{equation}
    \frac{\partial }{\partial t}H\leq -\frac{1}{n-1}H^2\label{Ricci>0}
\end{equation}
in the setting of non-negative Ricci curvature along normal flow $\frac{\partial }{\partial t}\Phi=\nu$. In non-negative Ricci curvature case, using this inequality, one can derive a Laplacian comparison theorem. Then, by applying techniques involving Busemann functions, it is possible to obtain a version of the Cheeger–Gromoll splitting theorem. This strategy can also be adapted to the static Einstein case via inequality (\ref{Riccati_equ}), although additional convexity assumptions on the conformal boundary are required to construct geodesic lines \cite{wang2024riccati}.

 Besides Busemann functions, another important technique for non-negative Ricci curvature is the instability of  minimal surface. Start from a closed minimal surface and run the normal flow $\frac{\partial }{\partial t}\Phi=\nu$. From (\ref{Ricci>0}) we know that $H$ remains non-positive and as a result area is non-increasing. So area-minimizing hypersurface does not exist unless metric splits.
 
Both the Busemann function technique and the theory of minimal surfaces are rooted in the fundamental principle that Ricci curvature controls mean curvature. Given that the Busemann function method can be successfully applied to the static Einstein setting, it is natural to ask whether minimal surface techniques can also be extended to this context. The answer is yes, and a lot results has been obtained \cite{Am,Ga,wang2024riccati}. However, all these results work on bounded minimal surfaces. In conformal compact spaces, there are a lot minimal surfaces which are not compact and touch conformal boundary orthogonally. it was shown in \cite{anderson1982complete,anderson1983complete} that in the hyperbolic space $(\mathbb{H}^n,g)$ with conformal boundary $(\mathbb{S}^{n-1},h)$, for any $k$ dimensional submanifold $\Gamma\subset \mathbb{S}^{n-1}$ with $k<n-1$, there exists a locally area-minimizing currents $F$ so that $F$ meets $\mathbb{S}^{n-1}$ along $\Gamma$ orthogonally. For these minimal surfaces, the area is infinite, so the methods mentioned above cannot be applied directly. To overcome this, the author introduces the concept of renormalized area, which has been intensively studied for PE manifolds,  and studies its instability, aiming to derive rigidity results.\\

Given a triple $(M^3,g,V)$ satisfying \textbf{Condition C}. According to \cite{graham1999conformal,Wo}, there exists a boundary defining function $x$ with $\lim xV=1$ such that near conformal boundary $g$ has expansion
\begin{equation}
    g=\frac{1}{x^2}(dx^2+h-\frac{R_{\Sigma}}{4}hx^2+h_3x^3+O(x^4))\label{expansion_x}
\end{equation}
In this expansion, $h$ is called Dirichlet data and the second order term $\frac{R_{h}}{4}hx^2$ is uniquely determined by $h$ from (\ref{static_Einstein_equ}). $h_3$ can not be calculated from $h$ from local information and is called Neumann data.

Denote $\bar{g}'=x^2g$. Let $F$ be a smooth surface such that $F\cap \Sigma=\Gamma$ for a smooth closed curve $\Gamma$ and $F\perp \Sigma$ in $\bar{g}'$ along $\Gamma$. Denote $F_{\epsilon}=F\cap \{p\in M:x(p)>\epsilon\}$, and compute its area, we have
\begin{equation}
    Area(F_{\epsilon},g)=\frac{Length(\Gamma,h)}{\epsilon}+const+O(\epsilon)\label{def_RenA_1}
\end{equation}
The const is defined to be the renormalized area of $F$ denoted by $RenA(F)$. According to \cite{graham1999volume}, $RenA$ is independant of the choice of boundary defining function and $h$ in the conformal class $[h]$.
\begin{theorem}\label{thm1}
    Let $(M,g,V)$ satisfying \textbf{Condition C}. Let $F_0$ be a minimal surface so that $F_0\cap \Sigma= \Gamma$ and $F_0\perp \Sigma$ (assume $\Gamma$ divides $\Sigma$ into two connected components). Let $F_t$ be the flow from $F_0$ satisfying (\ref{flow}) and $F_t\cap \Sigma =\Gamma_t$. Then
    \begin{equation}
        \begin{aligned}
            \frac{\partial }{\partial t}RenA(F(t))&=-\int_{\Gamma_0}3u_3d\theta_h\\
            \frac{\partial^2 }{\partial t^2}RenA(F(t))&=-\int_{F_0}V|b|^2dv_g+\int_{\Gamma_0}-3\kappa u_3+\frac{3}{4}(tr_{h}h_3+h_3(\partial s,\partial s))d\theta_h
         \end{aligned}
    \end{equation}
    where $u_3$ is the coefficient of $x^3$
in the expansion for the function $u$ which gives a graph
parametrization of $F_0$ over the vertical cylinder $\Gamma_0\times [0,x_0)$, $\kappa$ is the geodesic curvature of $\Gamma_0$ in $(\Sigma,h)$ and $\partial s$ is the normal vector to $\Gamma_0$ in $(\Sigma,h)$ and $b$ is the second fundamental form.
\end{theorem}
These quantities will be made more clear in section 3.

Based on the second variation formula, we are going to improve a rigidity result in \cite{GSW} for a warped product torus conformal boundary.

\begin{theorem}\label{thm2}
    Let $(M^3,g,V)$ satisfying \textbf{Condition C}. Suppose its compactification has topology $\bar{\mathbb{B}}^2\times S_s^1$, and $(\Sigma,h)$ is a torus given by $(\Sigma =S^1_\theta\times S^1_{s},h=f^2d\theta^2+ds^2)$ where $S_{\theta}=\partial \bar{\mathbb{B}}^2$. Then the Neumann data satisfies
    \begin{equation}
     \int_{\Sigma}tr_{h}h_3+h_3(\partial s,\partial s)dA_h\geq 0\notag
    \end{equation}    
    Besides, if 
    \begin{equation}
        \int_{\Gamma_{s_0}}tr_{h}h_3+h_3(\partial s,\partial s)d\theta\leq 0 \text{ for each }s_0\in S_{s}
    \end{equation}
    where $\Gamma_{s_0}=\{(\theta,s_0)\in \Sigma\}$ and $f$ is constant, then $(M,g,V)$ must be the Horowitz-Myers soliton given by (\ref{Horowitz-Myers}).
\end{theorem}

Existence of minimal surface is required to prove the theorem above.
\begin{definition}
    Let $(M^3,g)$ satisfy one of \textbf{Definition \ref{def_APE}}. Let $\Gamma$ a smooth closed curve in $\Sigma$. Define
    \begin{equation}
    \begin{aligned}
        C(\Gamma)&\coloneqq \{F: F \text{ is a complete surface, } F\cap \Sigma=\Gamma \text{ and }F\perp \Sigma\}\\
        i(\Gamma)&\coloneqq \inf\{RenA(F):F\in C(\Gamma)\}\label{def_admissible}
    \end{aligned}
    \end{equation}
\end{definition}
\begin{theorem}\label{thm3}
    Let $(M,g,V)$ be as in \textbf{Theorem \ref{thm2}}. Let $\Gamma_{s_0}=\{(\theta,s_0):\theta\in S_{\theta},s=s_0\}$. Then for each $s\in S_s$, there exists a minimal surface $F\in C(\Gamma(s_0))$ so that $RenA(F)=i(\Gamma(s))$.
\end{theorem}

This article is structured as follows: In Section 2, we present the main idea behind the calculation; in Section 3, we provide the necessary preliminaries; in Section 4, we compute the first and second variations of the renormalized area and prove Theorem 1; and finally, in Section 5, we prove the rigidity result. Unless stated otherwise $M$ will be 3-dimensional.

\section{Idea}
In this section we show the idea of computation and provide with a rough proof for Theorem 1.

\begin{figure}[htbp]
  \centering
  \begin{minipage}[b]{0.45\textwidth}
    \includegraphics[width=\textwidth]{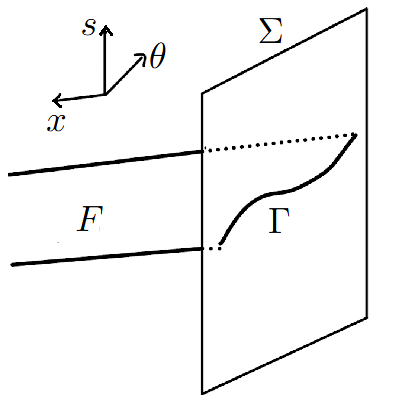}
    \caption{}
    \label{renA1}
  \end{minipage}
  \hfill
  \begin{minipage}[b]{0.45\textwidth}
    \includegraphics[width=\textwidth]{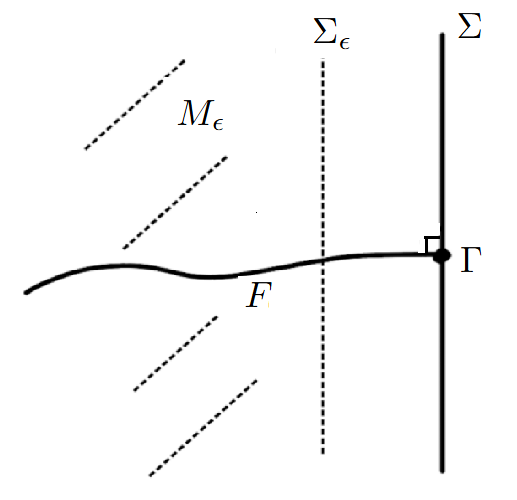}
    \caption{}
    \label{renA2}
  \end{minipage}
\end{figure}

Given $(M,g,V)$ and $F_t$ be as in Theorem \ref{thm1}. Let $r$ be a boundary defining function, $M_{\epsilon}=\{p\in M: r(x)>\epsilon\}$, $F_{t,\epsilon}=F_t\cap M_{\epsilon}$ and $\Gamma_t=F_t\cap \Sigma$. Recall that renormalized area is defined as
\begin{equation}
    RenA(F_t)=\lim_{\epsilon\rightarrow 0} \big(A(F_{t,\epsilon})-\frac{L(\Gamma_t)}{\epsilon}\big)\notag
\end{equation}
where $A(F_{t,\epsilon})=Area(F_{t,\epsilon},g)$ and $L(\Gamma)=Length(\Gamma,h)$. See figure \ref{renA1} and figure \ref{renA2}.
A straightforward idea is that
\begin{equation}
    \frac{\partial^2}{\partial t^2}RenA(F_t)=\lim_{\epsilon\rightarrow 0} \frac{\partial^2}{\partial t^2}\big(A(F_{t,\epsilon})-\frac{L(\Gamma_t)}{\epsilon}\big)\notag
\end{equation}
We set aside the issue of exchanging the order of limits for now. In order to compute the right hand side, we need to estimate
\begin{equation}
    (A(F_{t,\epsilon})-A(F_{0,\epsilon}))-\frac{1}{\epsilon}(L(\Gamma_t)-L(\Gamma_0))\label{2_subtrac}
\end{equation}

To estimate the subtraction in the first bracket, we present Figure \ref{idea1}:
\begin{figure}[htbp]
  \centering
  \begin{minipage}[b]{0.45\textwidth}
    \includegraphics[width=\textwidth]{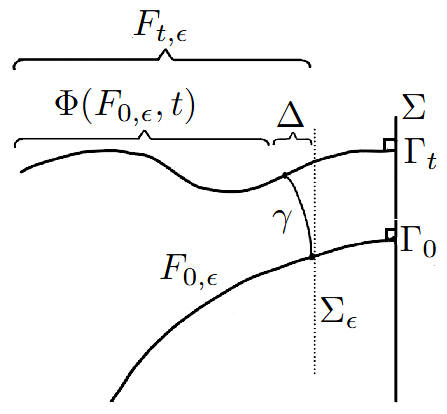}
    \caption{}
    \label{idea1}
  \end{minipage}
\end{figure}

\begin{figure}[htbp]
  \centering
  \begin{minipage}[b]{0.45\textwidth}
    \includegraphics[width=\textwidth]{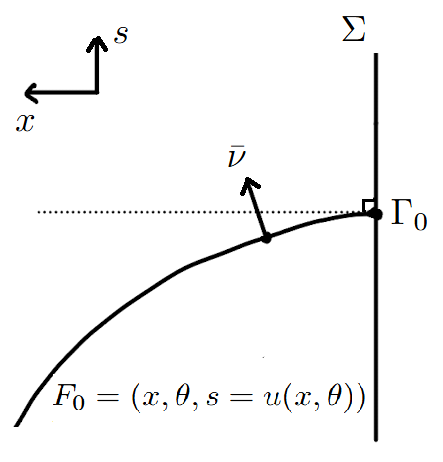}
    \caption{}
    \label{renA3}
  \end{minipage}
  \hfill
  \begin{minipage}[b]{0.45\textwidth}
    \includegraphics[width=\textwidth]{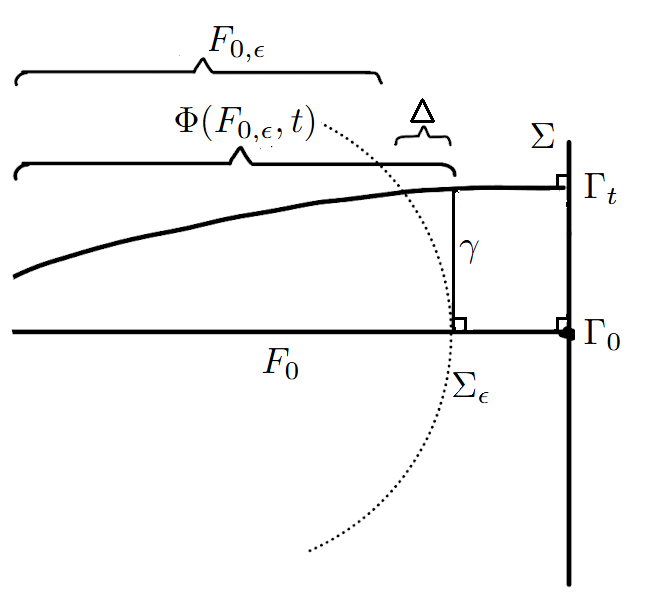}
    \caption{}
    \label{renA4}
  \end{minipage}
\end{figure}

$\gamma$ represents the trajectory of a point under the flow $\Phi$ and $\Phi(\cdot,t)$ maps $F_{0,\epsilon}$ to $\Phi(F_{0,\epsilon},t)$. Since we assume $F_0$ is a minimal surface, by (\ref{Riccati_equ}) $H(F_t) \leq 0$ and thus area is non-increasing along $\Phi$, so we have
\begin{equation}
    A(F_{t,\epsilon})-A(F_{0,\epsilon})=A(\Phi(F_{0,\epsilon},t))+A(\Delta)-A(F_{0,\epsilon})\leq A(\Delta)
\end{equation}
So key is to estimate $A(\Delta)$. Here the $A(\Delta)$ is signed area, and the sign is determined by the direction of $\gamma$: if $\gamma$ goes to the left of $\Sigma_{\epsilon}$, as in Figure \ref{renA3}, then $A(\Delta)>0$; and $A(\Delta)<0$ if $\gamma$ goes the other way. One observes that under conformal change the flow (\ref{flow}) becomes normal flow in $\bar{g}$, i.e. $\frac{\partial}{\partial t}\Phi=\bar{\nu}$ where $\bar{\nu}$ denotes normal vector in $\bar{g}=\frac{1}{V^2}g$. So we can use Jacobi equation to estimate the derivation of $\gamma$ from $\{p: r(p)=\epsilon\}$, thus estimating $A(\Delta)$.

We construct local coordinates as illustrated in Figure~\ref{renA1}. Roughly speaking, $\partial_{\theta}$ is tangent to $\Gamma$, $\partial_s$ lies along $\Sigma$ and is normal to $\partial_{\theta}$ while $\partial_x$ is normal to the conformal boundary $\Sigma$. In these coordinates, the rescaled metric $\frac{1}{V^2}g$ admits the following asymptotic expansion:
\begin{equation}
        \frac{1}{V^2}g\sim dx^2+h-\frac{R_{h}}{2}x^2h+(h_3+(tr_{h}h_3)h)x^3.\label{2_expansion_g}
    \end{equation}

    (The justification for this expansion will be provided in Section 3.)
In this coordinate system, the minimal surface $F$ can be represented as a graph over $(x,\theta)$, that is, $F=\{(x,\theta,u(x,\theta))\}$
as depicted in Figure~\ref{renA3}. Using the minimal surface equation $H=0$, the function $u$ also admits an asymptotic expansion:
\begin{equation}
    u\sim -\frac{\kappa}{2} x^2+u_3 x^3\label{2_expansion_u}
\end{equation}
where $\kappa$ denotes the geodesic curvature of $\Gamma$ in $\Sigma$.

It is important to note that in both expansions above, the coefficients $h_3$ and $u_3$	
  cannot be determined locally from the static Einstein equation or the minimal surface equation. These terms will play a crucial role in estimating the variation of the $x$-coordinate along a geodesic $\gamma$, denoted $\Delta \gamma^x$, and subsequently in estimating $A(\Delta)$.\\

The first variation in $t$ can be computed from the normal vector. Specifically, we have $\dot{\gamma}(0)=\bar{\nu}\sim (\kappa \epsilon+3u_3\epsilon^2,0,1)$, where this expression follows from the asymptotic expansions in equations~(\ref{2_expansion_g}) and~(\ref{2_expansion_u}). Let $\gamma^i$ denote the 
$i$-th coordinate function along the geodesic $\gamma$. In particular, we are interested in the 
$x$-component, which is given by
\begin{equation}
    \dot{\gamma}^x(0)\sim \kappa \epsilon+3u_3\epsilon^2\label{2_dx1}
\end{equation}

And the second variation in $t$ can be computed via Jacobi equation for geodesic  $\ddot{\gamma}^x=-\dot{\gamma}^i\dot{\gamma}^j\Gamma^x_{ij}$ and Einstein notation is used. Since $\dot{\gamma}(0)\sim (\kappa \epsilon+3u_3\epsilon^2,0,1)$, we have
\begin{equation}
    \ddot{\gamma}^x\sim\Gamma^x_{ss}\sim -\frac{R_h}{2}\epsilon+\frac{3}{2}(h_{3,ss}+tr_{h}h_3)\epsilon^2\label{2_dx2}
\end{equation}

It follows from (\ref{2_dx1})(\ref{2_dx2}) that
\begin{equation}
    \Delta \gamma^x\sim (\kappa \epsilon+3u_3\epsilon^2)t+\big(-\frac{R_h}{4}\epsilon+\frac{3}{4}(h_{3,ss}+tr_{h}h_3)\epsilon^2\big)t^2
\end{equation}
Again $F_t$ can be expressed as a graph over $(x,\theta)$, so its determinant can be computed to be $\sim \frac{1}{x^2}$. So
\begin{equation}
\begin{aligned}
    A(\Delta)&\sim \int_{\Gamma} \frac{1}{\epsilon^2}\Delta\gamma^x \\
    &\sim \frac{1}{\epsilon}\int_{\Gamma}(\kappa t-\frac{R_h}{4}t^2)+\int_{\Gamma}3u_3 t+\frac{3}{4}(h_{3,ss}+tr_{h}h_3)t^2
\end{aligned}
\end{equation}

Now compute the second two terms in (\ref{2_subtrac}). Since $\Sigma$ is umbilical in $(\bar{M},\bar{g})$, when restricted to $\Gamma_0$ the flow becomes normal flow in $(\Sigma,h)$. As a result
\begin{equation}
    \frac{1}{\epsilon}(L(\Gamma_t)-L(\Gamma_0))\sim \frac{1}{\epsilon}\int_{\Gamma_0}\kappa t-\frac{R_h}{4}t^2
\end{equation}

We now observe that, after performing the subtraction in equation~(\ref{2_subtrac}), the $\frac{1}{\epsilon}$
  term cancels out. As a result, we are nearly in a position to prove \textbf{Theorem~\ref{thm1}}. A more refined computation will yield the precise second variation.

To summarize, the core of the calculation lies in estimating \(A(\Delta)\). Its divergent part cancels with the change in \(\frac{1}{\epsilon} L(\Gamma)\), while its bounded part contributes to the first and second variations of the renormalized area, up to the additional term \(-\int V |b|^2\).\\

  For \textbf{Theorem~\ref{thm2}}, let \(\Gamma_{s_0} = \{ (\theta, s_0) \in \Sigma \}\). Then the function \(i(\Gamma_{s})\) is periodic in \(s\). Under our assumption, one can show that \(i(\Gamma_{s_0})\) is concave, and hence it must be constant. This constancy leads directly to the rigidity conclusion.\\

Finally, we use a slightly modified illustration (Figure~\ref{renA4}) to demonstrate how $h_3$ influences the renormalized area, particularly in the case where $R_h=0$. To better highlight the effect of 
$h_3$, we assume $u_2=u_3=0$. Under this assumption, if $h_3<0$, then the level sets $\Sigma_{\epsilon}$ are convex (the conformal boundary $\Sigma$ is umbilical). As a result, the geodesic $\gamma$ deviates to the right of $\Sigma_{\epsilon}$	, leading to the inequality $A(F_{t,\epsilon})<A(\Phi(F_{0,\epsilon},t))\leq A(F_{0,\epsilon})$. Thus, the negativity of $h_3$ (i.e., the convexity of $\Sigma_{\epsilon}$) implies the negativity of the second variation. The sign of \(h_3\) dictates the correction needed: if \(h_3 > 0\), \(\gamma\) deviates to the left, we must compensate for a deficit ; if \(h_3 < 0\), it deviates to the right, we cancel the excess.

\section{Preliminaries}
In this section, we provide the preliminaries required for the computation. We will always assume that the triple
 $(M^3,g,V)$ satisfies \textbf{Condition C}.

\subsection{Expansion for $g$ Near Conformal Boundary}

In this part, we derive the expansion of $\frac{1}{V^2}g$. To this end, we divide the argument into three lemmas:

(i) construct a suitable boundary defining function $x$ and obtain the expansion of the metric in terms of $x$;

(ii) establish the relationship between $\frac{1}{V}$ and $x$ using (\ref{static_Einstein_equ});

(iii) derive the expansion of $\frac{1}{V^2}g$.

\begin{lemma}\label{3_g_preliminary}
    There exists a boundary defining function with $\lim xV=1$ such that
    \begin{equation}
    g=\frac{1}{x^2}\big(dx^2+h-\frac{R_h}{4}x^2h+x^3h_3+o(x^3)\big) \label{expansion_1}
\end{equation}
\end{lemma}
\begin{proof}
We might lift $(M,g,V)$ to a Poincar\'e-Einstein manifold via
\begin{equation}
(M^3\times S^1,\hat{g} = V^2 d\theta^2 + g) 
\end{equation}
and the conformal boundary is $(\Sigma\times S^1,\hat{h}=h+d\theta^2)$ with boundary defining function $\frac{1}{V}$. By \cite{Lee1994TheSO}, there exists a boundary defining function $x$ so that $\hat{g}$ has the following expansion 
\begin{equation}
     \hat{g}=\frac{1}{x^2}\big(dx^2+\hat{h}(x)\big)\notag
\end{equation}
and $\hat{h}(x)$'s are a family of smooth metrics on $\Sigma\times S^1$. $x$ was obtained by solving a non-characteristic first order PDE, and therefore is unique. Due to the symmetry of $\hat{g}$ in $\theta$ direction, $x$ descends to a well-defined function on $M$. Furthermore, the expansion for $\hat{h}$ can be computed from Einstein equation as in
(see \cite{graham1999conformal})
\begin{equation}
    \hat{g}=\frac{1}{x^2}\big(dx^2+\hat{h}+x^2\hat{h}_2+x^3\hat{h}_3+o(x^3)\big)
\end{equation}
where $\hat{h}_2$ is given by
\begin{equation}
    \hat{h}_2=-(Ric(\hat{h})-\frac{R(\hat{h})}{4}\hat{h})
\end{equation}
Direct computation shows that
\begin{equation}
    \begin{aligned}
        \hat{h}(X,Y)&=-\frac{R_h}{4}h(X,Y) \quad \text{ for X,Y along } \Sigma\\
        \hat{h}(\partial \theta,\partial \theta)&=\frac{R_h}{4}
    \end{aligned}
\end{equation}
Now (\ref{expansion_1}) follows.
\end{proof}
\begin{definition}
    In \textbf{Lemma \ref{3_g_preliminary}}, $h_3$ can not be locally determined from static Einstein equation, and is defined as the \textbf{Neumann data}. The mass aspect function $\nu$ and the Wang mass $m(g)$ are defined as
    \begin{equation}
    \begin{aligned}
        \mu&\coloneqq 3tr_{h}h_3\\
        m(g)&\coloneqq \int_{\Sigma}\mu dv_h\label{3-1-mass}
    \end{aligned}
    \end{equation}
\end{definition}

\begin{lemma}
Near $\Sigma$, $\frac{1}{V}$ has the following asymptotic expansion
    \begin{equation}
        \frac{1}{V}=x-\frac{R_h}{8}x^3+\frac{tr_{h_0}h_3}{2}x^4+o(x^4)\label{3-V-expansion}
    \end{equation}
\end{lemma}
\begin{proof}
 For $p\in \Sigma$, we might find local normal coordinates $\{\partial_i\}_{i=1,2}$ w.r.t $h=h(0)$. Then extend it to coordinates system in a collar neighborhood $[0,x_0)\times \Sigma$ to$\{\partial_i\}_{i=0}^2$ where $\partial_0=\partial x$.
    Note that 
    \begin{equation}
        \Delta V=\sum_{0\leq i,j\leq 2}g^{ij}\nabla V^2(\partial_i,\partial_j)=x^2\nabla^2 V(\partial x,\partial x)+\sum_{1\leq i,j\leq 2}g^{ij}\nabla V^2(\partial_i,\partial_j)
    \end{equation}
        Plugging $(x\partial x,x\partial x)$ to both sides of (\ref{static_Einstein_equ}), $\nabla^2(\partial_x,\partial_x)$ term cancels and we get
    \begin{equation}
        \sum_{1\leq i,j\leq 2}g^{ij}\nabla V^2(\partial_i,\partial_j)=-VRic(x\partial x,x\partial x)\label{3-1-middle}
    \end{equation}
    For the left hand side, since $\{\partial_i\}_{i=1,2}$ are normal coordinates w.r.t $h=h(0)$, from the expansion (\ref{expansion_1}), we have
       \begin{equation}
       \begin{aligned}
           \Gamma_{ij}^0&=\frac{1}{2}g^{00}g_{ij,0}=\frac{1}{x}\delta_{ij}-\frac{1}{2}x^2h_{3,ij}+O(x^3)\\
           \Gamma_{ij}^k&=O(x^3) \text{  for }k=1,2
           \end{aligned}
    \end{equation}
    
         Let $u=\frac{1}{V}=x+bx^2+cx^3+dx^4+o(x^4)$, and $u_i=\partial_i u$ for $i=0,1,2$. For $i,j=1,2$ we have
            \begin{equation}
            \begin{aligned}
       \nabla^2 V(\partial_i,\partial_j)&=\big(-\frac{1}{u^2}\nabla^2 u(\partial_i,\partial_j)+\frac{2u_iu_j}{u^3}\big)\\
       &=\frac{1}{u^2}\big(-u_{ij}+\sum_{0\leq k\leq 2}u_k\Gamma_{ij}^k+\frac{2u_iu_j}{u}\big)\\
       &=\frac{1}{u^2}\big((\frac{1}{x}\delta_{ij}-\frac{1}{2}x^2h_{3,ij})u_0-u_{ij}+\frac{2u_iu_j}{u}+O(x^3)\big)
       \end{aligned}
       \end{equation}
           Since $g^{ij}=x^2(\delta_{ij}+\frac{R_h}{4}x^2\delta_{ij}-x^3h_{3,ij}+O(x^4))$, the left hand side of (\ref{3-1-middle}) is given by
           \begin{equation}
              \sum_{1\leq i,j\leq 2}g^{ij}\nabla V^2(\partial_i,\partial_j)=\frac{x^2}{u^2}\big( (\frac{2}{x}+\frac{R_h}{2}x-\frac{3}{2}tr_{h_0}(h_3)x^2)u_0+g^{ij}(-u_{ij}+\frac{2u_iu_j}{u})+O(x^3)\big)\label{3-1-middle-2}
           \end{equation}
    The right hand side of (\ref{3-1-middle}) can be computed via Riccati equation. Let $H$ be the mean curvature for level sets $\{x=c\}$. From (\ref{expansion_1}),
    \begin{equation}
    \begin{aligned}
        H&=\frac{1}{2}\sum_{i,j=1,2}g^{ij}x\frac{\partial}{\partial x}g_{ij}\\
        &=\frac{1}{2}\sum_{i,j=1,2}x^2(\delta_{ij}+\frac{R_h}{4}x^2\delta_{ij}-x^3h_{3,ij}+O(x^4))(-\frac{2}{x^2}\delta_{ij}+xh_{3,ij}+O(x^2))\\
        &=-2-\frac{R_h}{2}x^2+\frac{3}{2}tr_{h_0}(h_3)x^3+O(x^4)
    \end{aligned}
    \end{equation}
    And $Ric(x\partial_x,x\partial_x)$ is given by
    \begin{equation}
        \begin{aligned}
            Ric(x\partial_x,x\partial_x)=-x\frac{\partial}{\partial x}H-\frac{H^2}{2}=-2-\frac{3}{2}tr_{h_0}h_3x^3+O(x^4)\label{3-1-middle3}
        \end{aligned}
    \end{equation}
    Plug (\ref{3-1-middle-2}) and (\ref{3-1-middle3}) into (\ref{3-1-middle}), we arrive at
    \begin{equation}
         (2+\frac{R_h}{2}x^2-\frac{3}{2}tr_{h_0}(h_3)x^3)u_0+xg^{ij}(-u_{ij}+\frac{2u_iu_j}{u})+O(x^4)=\frac{u}{x}(2+\frac{3}{2}tr_{h_0}h_3x^3+O(x^4))\label{3-1-middle4}
    \end{equation}
Recall that $u=\frac{1}{V}=x+bx^2+cx^3+dx^4+o(x^4)$. By comparing at the 1st order term on both sides, we get $b=0$. As a result $u_{ij},\frac{u_iu_j}{u}$ are of order $O(x^3)$. Now (\ref{3-1-middle4}) becomes
\begin{equation}
         (2+\frac{R_h}{2}x^2-\frac{3}{2}tr_{h_0}(h_3)x^3)u_0+O(x^4)=\frac{u}{x}(2+\frac{3}{2}tr_{h_0}h_3x^3+O(x^4))\notag
\end{equation}
By comparing the remaining coefficients on both sides, we get $c=-\frac{R_h}{8}$ and $d=\frac{tr_{h_0}h_3}{2}$, and our lemma follows.
\end{proof}

Combining the two lemmas above, we get
\begin{lemma}
Near conformal boundary $\Sigma$, $\bar{g}=\frac{1}{V^2}g$ has the following
    \begin{equation}
        \frac{1}{V^2}g=(1-\frac{R_h}{4}x^2+tr_{h_0}h_3x^3)dx^2+h-\frac{R_h}{2}x^2h+(h_3+(tr_{h}h_3)h)x^3+O(x^4).\label{3_g_expansion}
    \end{equation}
\end{lemma}
\begin{example}
    Using Poincar\'e ball model, the hyperbolic space becomes
    \begin{equation}
        g_{\mathbb{H}^n}=\frac{4}{(1-r^2)^2}(dr^2+g_{\mathbb{S}^{n-1}})\notag
    \end{equation}
    where $g_{\mathbb{S}^{n-1}}$ is the round metric on $\mathbb{S}^{n-1}$. Using the $r=\frac{2-x}{2+x}$, it becomes
    \begin{equation}
        g_{\mathbb{H}^n}=\frac{1}{x^2}(dx^2+(1-\frac{x^2}{4})^2g_{\mathbb{S}^{n-1}})\notag
    \end{equation}
    For $n=3$, the Neumann data and potential are
    \begin{equation}
    \begin{aligned}
    h_3&=0\\
        \frac{1}{V}=\frac{x}{1+x^2/4}&=x-\frac{1}{4}x^3+O(x^4)\notag
        \end{aligned}
    \end{equation}
\end{example}
\begin{example}
    The Horowitz-Myers soliton (\ref{Horowitz-Myers}) can be expressed in the form (\ref{3_g_preliminary}) as follows:
\begin{equation}
    g=\frac{1}{x^2}\bigg(dx^2+(1+\frac{x^n}{4})^{4/n}[(\frac{1-x^n/4}{1+x^n/4})^2d\theta^2+\sum_{i=1}^{n-2}ds^2]\bigg)
\end{equation}
with $V=r=\frac{1}{x}(1+\frac{x^n}{4})^{2/n}$. For $n=3$, the Dirichlet data and Neumann data can be read off as
\begin{equation}
\begin{aligned}
    \text{Dirichlet data}&=d\theta^2+ds^2\\
    \text{Neumann data}&=-\frac{2}{3}d\theta^2+\frac{1}{3}ds^2\\
    \frac{1}{3}\mu&=tr_{h_0}h_3=-\frac{1}{6}\\
    \frac{1}{V}&=x(1-\frac{1}{3}x^3+O(x^4))
\end{aligned}
\end{equation}
So we have that
\begin{equation}
    h_3+(tr_{h_0}h_3)h_0=-d\theta^2
\end{equation}
And which coincides with that computed from (\ref{Horowitz-Myers}) and thus verifies (\ref{3-V-expansion})(\ref{3_g_expansion}). 
\end{example}

\subsection{Coordinates System}
In this section we are going to construct coordinates system along a given closed smooth curve $\Gamma$ and study the behavior of a minimal surface, and it also works for APE manifolds.

Let $\Gamma(\theta)\subset \Sigma$ be a closed curve that divides $\Gamma$ into two connected components. Let $\vec{n}$ denote the unit normal vector along $\Gamma$, pointing into one of the components. Using the exponential map in $(\Sigma,h)$, we obtain a collar neighborhood around $\Gamma$ given by
\begin{equation}
    N=\{exp_{\theta}(s\vec{n}):-s_0\leq s\leq s_0\}\notag
\end{equation}
And $(\theta,s)$ gives coordinates for $N$. Along normal flow, we have
\begin{equation}
    \begin{aligned}
        \frac{\partial}{\partial s}d\theta&=\kappa d\theta\\
        \frac{\partial^2}{\partial s^2}d\theta&=(\frac{\partial}{\partial s}\kappa)d\theta+\kappa\frac{\partial}{\partial s}d\theta\\
        &=(-\frac{R_h}{2}-\kappa^2)d\theta+\kappa^2d\theta=-\frac{R_h}{2}d\theta\notag
    \end{aligned}
\end{equation}
So within $N$
\begin{equation}
    h=ds^2+(1+\kappa s-\frac{R_h}{4}s^2+O(s^3))^2d\theta^2\label{3-expansion_h}
\end{equation}
Together, $(x,\theta,s)$ gives coordinate near $\Gamma$. And from section 2.3 of \cite{alexakis2010renormalized}, we have

\begin{lemma}
    Let $F$ a surface in $\bar{M}$ and $F\cap \Sigma=\Gamma$. Suppose $F$ can be expressed as a function $u(x,\theta)$ in the local coordinates above in a neighborhood of $\Gamma$, and $F$ is a minimal surface for $g$, then
    \begin{equation}
        u=-\frac{\kappa}{2}x^2+u_3x^3+O(x^4)\label{3-u}
    \end{equation}
    where $u_3$ can not locally be determined from the minimal surface equation.
\end{lemma}
\begin{remark}
    In \cite{alexakis2010renormalized}, the coefficient of the $x^2$ term appears as $\frac{\kappa}{2}$, whereas in our setting it is $-\frac{\kappa}{2}$. This discrepancy arises from a difference in sign convention and coordinate notation. Importantly, the coefficient $u_3$ cannot be determined locally from the minimal surface equation alone. Once $u_3$ is prescribed, the higher-order terms in the expansion can, in principle, be solved recursively. However, these formal expansions do not necessarily correspond to a complete minimal surface.
\end{remark}

\section{1st and 2nd Varation for Renormalized Area}

In this section we are going to prove \textbf{Theorem \ref{thm1}}. For reader's convenience, we restate the theorem here:
\begin{theorem}
    Let $(M,g,V)$ satisfying \textbf{Condition C}. Let $F_0$ be a minimal surface so that $F_0\cap \Sigma= \Gamma_0$ and $F_0\perp \Sigma$ (assume $\Gamma$ divides $\Sigma$ into two connected components). Let $F_t$ be the flow from $F_0$ satisfying (\ref{flow}) and $F_t\cap \Sigma =\Gamma_t$. Then
    \begin{equation}
        \begin{aligned}
            \frac{\partial }{\partial t}RenA(F(t))&=-\int_{\Gamma_0}u_3d\theta_h\\
            \frac{\partial^2 }{\partial t^2}RenA(F(t))&=-\int_{F_0}V|b|^2dv_g+\int_{\Gamma_0}-\kappa u_3+tr_{h}h_3+h_3(\partial s,\partial s)d\theta_h\label{thm1_result}
        \end{aligned}
    \end{equation}
    and the notations are the same as in section 3.
\end{theorem}

\begin{proof}
Before proceeding, we need to verify that the renormalized area is well defined for \(F_t\) for each \(t\). Note that the flow \(\partial_t \Phi = V \nu\) becomes the normal flow with respect to \(\bar{g}\). From equation~\eqref{3_g_expansion}, we can glue two \((\bar{M}, \bar{g})\) along \(\Sigma\) to form a \(C^2\) manifold \((\bar{M}^*, \bar{g}^*)\). At the same time, the flow \(\Phi\) is lifted to a \(C^1\) normal flow \(\Phi^*\) on \(\bar{M}^*\). By symmetry, \(\Phi^*\) is perpendicular to \(\Sigma\), and thus the renormalized area is well-defined for it. Besides, the trajectory from $\Gamma_0$ stays within $\Sigma$.

Let $M_{\epsilon}=\{x>\epsilon\} $, $\Sigma_{\epsilon}=\partial M_{\epsilon} $, $F_{t}=\Phi(\cdot,t)$ and $F_{t,\epsilon}=F_t\cap M_{\epsilon}$ (see Figure \ref{idea1}). Fix a small $\epsilon_0$, and we have
\begin{equation}
    RenA(F_t)=A(F_{t,\epsilon_0})-\frac{L(\Gamma_t)}{\epsilon_0}+\lim_{\epsilon\rightarrow 0}\big((A(F_{t,\epsilon}\setminus F_{t,\epsilon_0})-(\frac{1}{\epsilon}-\frac{1}{\epsilon_0})L(\Gamma_t)\big)\notag
\end{equation}
where $A(F_{t,\epsilon_0})=Area(F_{t,\epsilon_0},g)$ and $L(\Gamma_t)=Length(\Gamma_t,h)$. For simplicity, the dependence on $g$ and $h$ will be dropped since there will be no confusion. 
Using the strategy in Section 2, i.e. estimating $A(\Delta)\coloneqq A(F_{t,\epsilon_0})-A(\Phi(F_{0,\epsilon},t))$, $RenA(F_t)$ can be rewritten as

\begin{equation}
\begin{aligned}
    RenA(F_t)=&A(\Phi(F_{0,\epsilon},t))\!+\!\bigg(\!A(\Delta)\!-\!  \frac{L(\Gamma_t)}{\epsilon_0}\!\bigg)\!+\!\lim_{\epsilon\rightarrow 0}\big((A(F_{t,\epsilon}\!\setminus\! F_{t,\epsilon_0})\!-\!(\frac{1}{\epsilon}\!-\!\frac{1}{\epsilon_0})L(\Gamma_t)\big)\\
    =&\quad I\qquad \quad +  \qquad \quad  II\qquad \qquad +\qquad \qquad \qquad III  \label{4_initial}
\end{aligned}
\end{equation}
we are going to show that 
\begin{equation}
    \begin{aligned}
        \frac{\partial^2}{\partial t^2}I&= -\int_{F_0}V|b|^2dv_g\\
        \frac{\partial^2}{\partial t^2}II&=\int_{\Gamma_0}-3\kappa u_3+\frac{3}{4}(tr_{h}h_3+h_3(\partial s,\partial s))d\theta_h+O(\epsilon_0)\\
        \frac{\partial^2}{\partial t^2}III&=O(\epsilon_0)
    \end{aligned}
\end{equation}
Let $\epsilon_0\rightarrow 0$ and the theorem follows.

The major difficulty is to estimate  $\frac{\partial^2}{\partial t^2}II$. To do this, the proof is divided into 3 parts: i) evaluating the behavior of the flow $\partial_t \Phi=\bar{\nu}$ near conformal boundary; ii) compute $\frac{\partial^2}{\partial t^2}II$; iii) compute $\frac{\partial^2}{\partial t^2}I$ and compute $\frac{\partial^2}{\partial t^2}III$.

Following the coordinate in section 3.2. Let $\Sigma_{\epsilon}=\{x=\epsilon\}$ and $\Gamma_{t,\epsilon}=F_t\cap \Sigma_{\epsilon}$. Let $\gamma(t)$ a geodesic for $\bar{g}$ emanating from $F_0$ and normal to $F_0$. Then $A(F_{t,\epsilon_0})-A(\Phi(F_{0,\epsilon},t))$ is the area of the surface that lies between $\Phi(\Gamma_{0,\epsilon_0},t)$ and $\Gamma_{t,\epsilon_0}$, denoted by $\Delta$, with a sign determined by the following: if $x>\epsilon_0$ in this region ($\gamma$ deviated to the left of $\Sigma_{\epsilon}$, see Figure~\ref{idea1}), then the area is taken to be positive; otherwise, it is negative.

View $F_t$ as a graph over $(x,\theta)$ with area element $\text{Det}_{t}$, and projects $\Delta$ to a region $\Delta_0$ in the $(x,\theta)$ plane, the above can be expressed as
\begin{equation}
    A(F_{t,\epsilon_0})-A(\Phi(F_{0,\epsilon},t))=\int_{\Delta_0}\text{sgn}(x-\epsilon_0)\text{Det}_tdxd\theta\notag
\end{equation}
Further, the projection of $\Phi(\Gamma_{0,\epsilon},t)$ to $(x,\theta)$ plane can be expressed as a graph over $\theta$, i.e. $(w(\theta),\theta)$. The above equation becomes 
\begin{equation}
    A(F_{t,\epsilon_0})-A(\Phi(F_{0,\epsilon},t))=\int_{\Gamma}\int_{\epsilon_0}^{w(\theta)}\text{Det}_tdxd\theta\notag
\end{equation}

   \noindent \textit{i):Estimate $\partial_t \Phi=\bar{\nu}$ near conformal boundary}

   Let $\Gamma_{\epsilon}=F\cap \Sigma_{\epsilon}$, we need to find the behavior of geodesics emanating from $\Gamma_{\epsilon}$ and normal to $F$. Let $\gamma$ be such a geodesic. $F$ can locally be expressed as graph $(x,\theta,s=u(x,\theta))$. We can choose
    \begin{equation}
        \partial_1=(1,0,\frac{\partial u}{\partial x})=-\kappa x+3u_3x^2+O(x^3);\quad \partial_2=(0,1,\frac{\partial u}{\partial \theta})=O(x^2).
    \end{equation}
    A direct algebra computation shows that one vector perpendicular to $F$ in $\bar{g}$ is given by $(a,b,c)$ where
   \begin{equation}
       \begin{aligned}
           a&=\bar{g}_{x\theta}\bar{g}_{\theta s}-\bar{g}_{\theta\theta}\bar{g}_{x s}+\frac{\partial u}{\partial x}(\bar{g}^2_{\theta s}-\bar{g}_{\theta\theta}\bar{g}_{s s})+\frac{\partial u}{\partial \theta}(\bar{g}_{x\theta}\bar{g}_{s s}-\bar{g}_{xs}\bar{g}_{\theta s})\\
           b&=\bar{g}_{x\theta}\bar{g}_{x s}-\bar{g}_{xx}\bar{g}_{\theta s}+\frac{\partial u}{\partial x}(\bar{g}_{x\theta}\bar{g}_{s s}-\bar{g}_{xs}\bar{g}_{\theta s})+\frac{\partial u}{\partial \theta}(\bar{g}^2_{xs}-\bar{g}_{xx}\bar{g}_{s s})\\
           c&=\bar{g}_{xx}\bar{g}_{\theta\theta}-\bar{g}^2_{x\theta}+\frac{\partial u}{\partial x}(\bar{g}_{xs}\bar{g}_{\theta\theta}-\bar{g}_{x\theta}\bar{g}_{\theta s})+\frac{\partial u}{\partial \theta}(\bar{g}_{xx}\bar{g}_{\theta s}-\bar{g}_{x\theta}\bar{g}_{x s})\\
       \end{aligned}
   \end{equation}
   From (\ref{3_g_expansion}) and (\ref{3-u}), we have
   \begin{equation}
       \begin{aligned}
           a&=-\frac{\partial u}{\partial x}\bar{g}_{\theta\theta}\bar{g}_{s s}+O(x^3)=\kappa x-3u_3x^2+O(x^3)\\
           b&=-\frac{\partial u}{\partial \theta}\bar{g}_{xx}\bar{g}_{s s}+O(x^3)=O(x^2)\\
           c&=\bar{g}_{xx}\bar{g}_{\theta\theta}+O(x^3)=1+O(x^2)
       \end{aligned}
   \end{equation}
   After uniformalization, the unit normal vector w.r.t. $\bar{g}$ is given by
   \begin{equation}
       \dot{\gamma}(0)=\bar{\nu}=(\kappa x-3u_3x^2+O(x^3),O(x^2),1+O(x^2))\label{4-normal}
   \end{equation}
This shows the coefficients of $\Delta \gamma^x$ for $t$ is given by $\kappa x-3u_3x^2+O(x^3))$. To compute coefficients for $t^2$, we need Jacobi equation $\ddot{\gamma}^k=-\dot{\gamma}^i\dot{\gamma}^j\Gamma^k_{ij}$ ($\Gamma^k_{ij} $ denotes the Christoffel symbols for $\bar{g}$). From (\ref{3_g_expansion}) and (\ref{4-normal}), direct computation shows that $\ddot{\gamma}^x(0)=-\Gamma_{ss}^{x}+O(x^3)=-\frac{R_h}{2}x+\frac{3}{2}(h_{3,ss}+tr_{h}h_3)x^2+O(x^3),
    \ddot{\gamma}^{\theta}(0)=O(x^2)$ and $
    \ddot{\gamma}^s(0)=O(x^2)$.

(One way to see this is that (\ref{3_g_expansion}) can be rewritten as $\bar{g}=dx^2+ds^2+(1+2\kappa s+O(s^2))d\theta^2+O(x^2)$. If $i,j,k\neq x$, then $\partial_ig_{jk}=O(x^3)$ except $\partial_s g_{\theta\theta}=2\kappa$; if exactly one of the indices $i,j,k$ is equal to $x$, then $\partial_ig_{jk}=O(x)$; if exactly two of indices of $i,j,k$ are equal to $x$, then $\partial_ig_{jk}=O(x^2)$.)

It follows that
\begin{equation}
    \begin{aligned}
        \gamma^{x}(t)\!&=\!x \!+\! (\kappa x \!-\! 3u_3x^2 \!+\! O(x^3))t \!+\! \left(-\frac{R_h}{4}x \!+\! \frac{3}{4}(h_{3,ss} \!+\! \mathrm{tr}_{h}h_3)x^2 \!+\! O(x^3)\right)t^2 \!+\! O(t^3)\\
        \gamma^{\theta}(t)\!&=\theta+\!O(x^2)t \!+\! O(x^2)t^2 \!+\! O(t^3)\\
        \gamma^s(t)\!&= \!\! O(x^2) \!+\! ( 1 \!+\! O(x^2))t \!+\! O(x^2) t^2 \!+\! O(t^3)\label{estimate_gamma}
    \end{aligned}
\end{equation}
From this, we can estimate $w(\theta)$:
\begin{lemma}\label{lemma_w}
    \begin{equation}
        w(\theta)=\!\epsilon_0 \!+\! (\kappa \epsilon_0 \!-\! 3u_3\epsilon_0^2 \!+\! O(\epsilon_0^3))t \!+\! \left(-\frac{R_h}{4}\epsilon_0 \!+\! \frac{3}{4}(h_{3,ss} \!+\! \mathrm{tr}_{h}h_3)\epsilon_0^2 \!+\! O(\epsilon_0^3)\right)t^2 \!+\! O(t^3)\label{4-w(theta)}
    \end{equation}
\end{lemma}
\begin{proof}


The function $\gamma^x$ does not directly yield $w(\theta)$, as we must account for the sliding in the $\theta$-direction during the flow. To make this precise, consider the composition of the flow map $\Phi(\cdot, t)$ with the projection onto the $(x, \theta)$-plane, applied to $\Gamma_{0,\epsilon_0}$. This yields a map
$(\epsilon_0, \theta) \mapsto (\tilde{x}(\theta), \tilde{\theta}(\theta))$,
where $\tilde{x}(\theta)$ and $\tilde{\theta}(\theta)$ are given by $\gamma^x(t)$ and $\gamma^\theta(t)$, respectively, as in equation~\eqref{estimate_gamma}, with $x$ evaluated at $\epsilon_0$.

Since $\tilde{\theta}$ depends smoothly on $\theta$, this defines a diffeomorphism from $\theta$ to $\tilde{\theta}$. By definition, we have $w(\tilde{\theta}) = \tilde{x}$. However, since $\tilde{x}$ is expressed as a function of $\theta$, we must invert the diffeomorphism to express $w$ purely in terms of $\tilde{\theta}$. Using the expansion
$\tilde{\theta} = \theta + O(\epsilon_0^2)t + O(t^2),$
we obtain the estimate
\[
\kappa(\tilde{\theta}) + C \epsilon_0^2 t \leq \kappa(\theta) \leq \kappa(\tilde{\theta}) + C \epsilon_0^2 t,
\]
 The remaining terms in $\gamma^x(t)$ can be estimated in the same way, and the lemma follows.

\end{proof}
\begin{remark}
Actually, from the computation above, we are lucky enough that the shift in $\theta$ does not affect the lower-order terms. However, this change of variable is still necessary for accuracy in the expansion.

\end{remark}

Next, we need to estimate $\text{Det}_t$.
\begin{lemma}\label{lemma_det_expansion}
$\text{Det}_t$ has the following asymptotica expansion
    \begin{equation}
        \text{Det}_t=\frac{1}{x^2}\big(1+O(x^2)+(\kappa+O(x^2))t+(-\frac{R_h}{2}+O(x^2))t^2+O(t^3)\big)
    \end{equation}
\end{lemma}
\begin{proof}
    We need to express $F_t$ as a graph over $(x,\theta)$ by using (\ref{estimate_gamma}). Just as in \textbf{Lemma \ref{lemma_w}}, in the (\ref{estimate_gamma}) $\gamma^s$ is expressed as a function of $(x,\theta)$, but at the point $(\gamma^x,\gamma^{\theta},\gamma^x)$ the right coordinate is $(\gamma^x,\gamma^{\theta})$. Locally there is a diffeomorphism $(\gamma^x,\gamma^{\theta})=\phi(x,\theta)$ via the first two equations in (\ref{estimate_gamma}), and $F_t$ can be expressed as a graph over $(x,\theta)$ via $u_t$ where $u_t(\gamma^x,\gamma^{\theta})=\gamma^s(\phi^{-1}(\gamma^x,\gamma^{\theta}))$. For simplicity of the formula, we use the form $u_t(x,\theta)$ instead of $u_t(\gamma^x,\gamma^{\theta})$
    \begin{equation}
        u_t(x,\theta)=O(x^2) + ( 1 + O(x^2))t + O(x^2) t^2 + O(t^3)
    \end{equation}
    So we have
    \begin{equation}
    \begin{aligned}
    \partial_1&=(1,0,\partial_xu_t)=\big(1,0,x(O(1)+O(1)t+O(1)t^2+O(t^3))\big)\\
    \partial_2&=(0,1,\partial_{\theta}u_t)=\big(0,1,x^2(O(1)+O(1)t+O(1)t^2+O(t^3))\big)
    \end{aligned}
    \end{equation}
    Using (\ref{3_g_expansion})
    \begin{equation}
        \begin{aligned}           g(\partial_1,\partial_1)&=\frac{1}{x^2(}1+O(x^2)+O(x^2)t+O(x^2)t^2)\\
        g(\partial_2,\partial_2)&=\frac{1}{x^2}(1+O(x^2)+(2\kappa+O(x^2))t+(-\frac{R_h}{2}+O(x^2))t^2+O(t^3)))\\
        g(\partial_1,\partial_2)&=O(1)(1+t+t^2+O(t^3))
        \end{aligned}
    \end{equation}
Recall that $\text{Det}_t=\sqrt{g(\partial_1,\partial_1)g(\partial_2\partial_2)-g(\partial_1,\partial_2)^2}$, and the lemma follows.
\end{proof}

\noindent\textit{ii) First and second derivatives for $II$}

Now we are in the position to estimate $A(F_{t,\epsilon_0})-A(\Phi(F_{0,\epsilon},t))$
\begin{equation}
    \begin{aligned}
        &\quad\int_{\Gamma_{0,\epsilon_0}}\int_{\epsilon_0}^{w(\theta)}\text{Det}_tdxd\theta\\
        &=\int_{\Gamma_{0,\epsilon_0}}\int_{\epsilon_0}^{w}\frac{1}{x^2}\bigg((1+O(x^2)+(\kappa+O(x^2))t+(-\frac{R_h}{2}+O(x^2))t^2+O(t^3)\bigg)dxd\theta\\        &=\int_{\Gamma_{0,\epsilon_0}}\int_{\epsilon_0}^{w}\frac{1}{x^2}\big(1+\kappa t-\frac{R_h}{2}t^2\big)+O(1)(t+t^2)+O(t^3)dxd\theta\\
        &=\int_{\Gamma_{0,\epsilon_0}}(1+\kappa t-\frac{R_h}{2}t^2)\big(\frac{1}{\epsilon_0}-\frac{1}{w(\theta)}\big)+O(\epsilon_0)(t+t^2)+O(t^3)\label{4_middle}
    \end{aligned}
\end{equation}
Using (\ref{4-w(theta)}), we have
\begin{equation}
    \begin{aligned}
        &\quad \frac{1}{\epsilon_0}-\frac{1}{w(\theta)}\\
        &=\frac{1}{\epsilon_0} \bigg( 1 \!-\! \frac{1}{1 \!+\! (\kappa \!-\! 3u_3\epsilon_0 \!+\! O(\epsilon_0^2))t \!+\! \left(-\frac{R_h}{4} \!+\! \frac{3}{4}(h_{3,ss} \!+\! tr_{h} h_3)\epsilon_0+O(\epsilon_0^2)\right)t^2 \! \!+\! O(t^3)} \bigg)\\
        &=\frac{1}{\epsilon_0}\bigg(\! (\kappa \!-\! 3u_3\epsilon_0 \!+\! O(\epsilon_0^2))t \!+\! \left(-\frac{R_h}{4}\!-\kappa^2 \!+\! \frac{3}{4}(h_{3,ss} \!+\! tr_{h}h_3)\epsilon_0 \!+\! O(\epsilon_0^2)\right)t^2 \!+\! O(t^3)\bigg)
    \end{aligned}
\end{equation}
where in the last equality we used $1-\frac{1}{1+\delta}=\delta-\delta^2+\cdots$. The $\delta^2$ term will contribute to the extra $-\kappa^2t^2$ term. Following (\ref{4_middle}), we have
\begin{equation}
    \begin{aligned}
&\quad\int_{\Gamma_{0,\epsilon}}\int_{\epsilon_0}^{w(\theta)}\text{Det}_tdxd\theta\\
&=\big(\int_{\Gamma_{0,\epsilon}}(\frac{\kappa}{\epsilon_0}-3u_3)d\theta+O(\epsilon_0)\big)t\\
&\quad+\big(\int_{\Gamma_{0,\epsilon}}(-\frac{R_h}{4\epsilon_0}-3u_3\kappa+\frac{3}{4}(h_{3,ss} \!+\! tr_{h}h_3))d\theta+O(\epsilon_0)\big)t^2+O(t^3)\label{4_middle2}
    \end{aligned}
\end{equation}
By (\ref{3-expansion_h}),
\begin{equation}
    L(\Gamma_{t})-L(\Gamma_{0})=\kappa t-\frac{R_h}{4}t^2+O(t^3)
\end{equation}
Combining this with equation~\eqref{4_middle2}, we compute $A(\Delta) - \frac{1}{\epsilon_0} \left( L(\Gamma_t) - L(\Gamma_0) \right)$,
and observe a cancellation for the blow up terms: the linear term $\frac{1}{\epsilon_0} \kappa t$ and the quadratic term $-\frac{R_h}{4\epsilon_0} t^2$. As a result, we arrive at

\begin{equation}
    \begin{aligned}
        \frac{\partial}{\partial t}II&=\int_{\Gamma_{0,\epsilon}}-3u_3d\theta+O(\epsilon_0)\\
        \frac{\partial^2}{\partial t^2}II&=\int_{\Gamma_{0,\epsilon}}(-3u_3\kappa+\frac{3}{4}(h_{3,ss} \!+\! tr_{h}h_3))d\theta+O(\epsilon_0)\label{4_dIIdt}
    \end{aligned}
\end{equation}

\noindent\textit{iii) First and second derivatives for $I$ and $III$}

Since $F_0$ is a minimal surface, 
\begin{equation}
    \frac{\partial}{\partial t}I=0\label{4_dIdt}
\end{equation}
Along $\partial_t\Phi=V\nu$, mean curvature satisfies the evolution equation
\begin{equation}
\begin{aligned}
    \frac{\partial}{\partial t}H=-\Delta_{F_t}V-V|b|^2-VRic_{g}(\nu,\nu)
\end{aligned}
\end{equation}
Proceed with static Einstein equation (\ref{static_Einstein_equ}),
\begin{equation}
    \begin{aligned}
        \frac{\partial}{\partial t}H&=-\big(\Delta_{g}V-\nabla^2_gV(\nu,\nu)-Hg(\nabla V,\nu)+V|b|^2+VRic_{g}(\nu,\nu)\big)\\
        &=-V|b|^2+Hg(\nabla V,\nu)\label{4-dHdt}
    \end{aligned}
\end{equation}
Since $F_0$ is a minimal surface, the second derivative is
    \begin{equation}
        \frac{\partial^2}{\partial t^2}I= \frac{\partial}{\partial t}\int_{F_{0,\epsilon_0}} VHdv_g= -\int_{F_{0,\epsilon_0}}V^2|b|^2dv_g\label{4-dIdtt}
\end{equation}

 To compute derivatives for $III$, we need the following
 \begin{equation}
     (A(F_{t,\epsilon}\setminus F_{t,\epsilon_0})-(\frac{1}{\epsilon}-\frac{1}{\epsilon_0})L(\Gamma_t)=\int_{S^1}\int_{\epsilon}^{\epsilon_0}(\text{Det}_t-\frac{L(\Gamma_t)}{x^2})dxd\theta\notag
 \end{equation}
 From \textbf{Lemma \ref{lemma_det_expansion}} and (\ref{3-expansion_h}), we see that again the cancellation of blow-up terms for $\text{Det}_t-\frac{L(\Gamma_t)}{x^2}=O(1)(t+t^2)+O(t^3)$, and thus 
 \begin{equation}
     \frac{\partial}{\partial t}III=\frac{\partial^2}{\partial t^2}III=O(\epsilon)\label{4_dIIIdt}
 \end{equation}

 From (\ref{4_dIIdt})(\ref{4_dIdt})(\ref{4-dIdtt})(\ref{4_dIIIdt}), letting $\epsilon_0\rightarrow 0$ in (\ref{4_initial}) and we are done.
\end{proof}

In \cite{alexakis2010renormalized}, it was shown that if $(M^n, g)$ is a Poincar\'e--Einstein manifold with boundary defining function (bdf) $x$ and conformal boundary $\Sigma$ (with $n$ not necessarily equal to $3$), and if $F^2$ is a complete surface with asymptotic boundary $\Gamma$, then
\begin{equation}
    \operatorname{RenA}(F) = -2\pi \chi(F) + \int_F \left( \frac{1}{4}|H|^2 - |\mathring{b}|^2 \right) dA + \int_F W_{1212}\, dA,
\end{equation}
where $\mathring{b}$ is the traceless second fundamental form and $W_{1212}$ is the Weyl curvature of $g$ evaluated on any orthonormal basis of $TF$.

Based on this, it was shown that under the flow $\partial_t \Phi = \frac{1}{x} \phi_{-1} + O(1)$, we have
\begin{equation}
    \frac{\partial}{\partial t} \operatorname{RenA}(\phi(F, t)) = \int_\Gamma -3 u_3 \Phi_{-1},\label{4_Alex}
\end{equation}
where $u_3$ is as defined earlier. Furthermore, the second variation formula is computed if $F_t$ are all minimal surfaces. These argument relies on the specific structure of the Poincar\'e--Einstein manifold. In our case, since we consider a normal flow, we have $\Phi_{-1} = 1$, and (\ref{4_Alex}) agrees with the statement of \textbf{Theorem~\ref{thm1}}. Actually, our approach extends to more general settings beyond the setting of \textbf{Theorem~\ref{thm1}} (see the theorem below), but the second variation formula is restricted to static Einstein manifold with flow $\partial_t\Phi=V\nu$.

\begin{theorem}\label{thm5}
    Let $(M^3,g)$ be APE, $\Gamma$ a closed smooth curve in $\Sigma$. Pick local coordinates as in section 3. Suppose $F_0\in C(\Gamma)$ and satisfy
    \begin{equation}
        u=-\frac{\kappa}{2}x^2+u_3x^3+O(x^4)\label{4-2_u}
    \end{equation}
    where $u$ is the graph function for $F$. Let $\partial_t\Phi(\cdot,t)=\phi\nu$ be a variation so that near the conformal boundary $\phi$ takes form
    \begin{equation}
        \phi=\frac{\phi_{-1}}{x}+O(x)\label{4-2-flow}
    \end{equation}
    where $\phi_{-1}$ is a bounded and smooth function on $\Gamma_0$, then we have
    \begin{equation}
    \frac{\partial}{\partial t} \operatorname{RenA}(\phi(F, t)) = \int_{F_0}\phi HdV_g+\int_\Gamma -3 u_3 \phi_{-1}d\theta\label{thm5_result}
\end{equation}
\end{theorem}
\begin{proof}
    The proof is almost the same, so we only outline the sketch. Using (\ref{4_initial}) we split $RenA(F_t)$ into three part. Since $(M^3,g)$ is APE, we can find a bdf $x$ such that
    \begin{equation}
        g=\frac{1}{x^2}(dx^2+h+x^2h_2+O(x^3))
    \end{equation}
    $h_2$ depends on $h_0$, it's not needed to compute the first variation of $RenA$, and the only thing we need is that there's no $xh_1$ term. 
    By our assumption on $u$, (\ref{4-normal}) still holds, and then we have
    \begin{equation}
        \begin{aligned}
            \gamma^{x}(t)&=x + (\kappa \phi_{-1}x - 3u_3\phi_{-1}x^2 + O(x^3))t \!+\! O(t^2)\\
        \gamma^{\theta}(t)&=\theta+\!O(x^2)t + O(t^2)\\
        \gamma^s(t)&=  O(x^2) + ( \phi_{-1} + O(x^2))t + O(t^2)\label{4-2-estimate_gamma}
        \end{aligned}
    \end{equation}
It follows that 
\begin{equation}
        w=\epsilon_0+(\kappa \phi_{-1}\epsilon_0 - 3u_3\phi_{-1}\epsilon_0^2 + O(\epsilon_0^3))t+O(t^2)
\end{equation}
    and 
        \begin{equation}
    \begin{aligned}
    \partial_1&=\big(1,0,x(O(1)+O(1)t+O(1)t^2+O(t^3))\big)\\
    \partial_2&=\big(0,1,x^2(O(1)+O(1)t+O(1)t^2+O(t^3))\big)\\
    \Longrightarrow Det_t&=\frac{1}{x^2}\bigg(1+O(x^2)+(\kappa+O(x^2))t+O(t^2)\bigg)\label{4-2-det}
    \end{aligned}
    \end{equation}
    And it follows that
    \begin{equation}
        \begin{aligned}
            A(F_{t,\epsilon_0})-A(\Phi(F_{0,\epsilon},t))&=\int_{\Gamma_{\epsilon_0}}\int_{\epsilon_0}^{w(\theta)}Det_tdxt\theta\\
            &=\int_{\Gamma_{\epsilon_0}}\big(\frac{1}{\epsilon_0}\phi_{-1}\kappa-3u_3\phi_{-1}+O(\epsilon_0)\big)td\theta+O(t^2)
        \end{aligned}
    \end{equation}
    On $\Sigma$, we have $L(\Gamma_t)-L(\Gamma_0)=\int\phi_{-1}\kappa d\theta$, so we see the cancellation in the blow-up term again:
    \begin{equation}
        \begin{aligned}
            &A(F_{t,\epsilon_0})-A(\Phi(F_{0,\epsilon},t))-\frac{1}{\epsilon_0}(L(\Gamma_t)-L(\Gamma_0))=\int_{\Gamma_{\epsilon_0}}(-3u_3\phi_{-1}+O(\epsilon_0))td\theta+O(t^2)\\
            &\hspace{4cm}\Longrightarrow \frac{\partial}{\partial t}II=\int_{\Gamma_{\epsilon_0}}-3u_3\phi_{-1}d\theta+O(\epsilon_0)
        \end{aligned}
    \end{equation}
    $I$ and $III$ can also be estimated:
    \begin{equation}
    \begin{aligned}
        \frac{\partial}{\partial_t}I&=\int_{F_{0,\epsilon_0}}\phi Hdv_g\\
         \frac{\partial}{\partial_t}III&=\frac{\partial}{\partial_t}\lim_{\epsilon\rightarrow 0}\int_{\Gamma_t}\int_{\epsilon}^{\epsilon_0}\text{Det}_t-\frac{1}{x}L(\Gamma_t)dxd\theta\\
         &=\frac{\partial}{\partial_t}\lim_{\epsilon\rightarrow 0}\int_{\Gamma_t}\int_{\epsilon}^{\epsilon_0}O(1)+O(1)t+O(t^2)dxd\theta=O(\epsilon_0)\label{4-2-tail}
    \end{aligned}
    \end{equation}
    And the proof is finished.
\end{proof}

\begin{remark}
i) We did not assume that \(F_0\) is a minimal surface, but we still impose the condition (\ref{4-2_u}). This assumption is necessary in order to ensure the cancellation with the blow-up term \(\frac{1}{\epsilon_0} (L(\Gamma_t) - L(\Gamma_0))\).

ii)    In general, $\phi=\frac{\phi_{-1}}{x}+\phi_0+O(x)$. But the theorem above is restriced to the case $\phi_0=0$. If $\phi_0$ is not $0$, then (\ref{4-2-estimate_gamma})(\ref{4-2-det}) transform sequentially into
    \begin{equation}
        \begin{aligned}
            \gamma^s(t)&=  O(x^2) + ( \phi_{-1} +\phi_0x+ O(x^2))t + O(t^2)\\
            \Longrightarrow \partial_1&=\big(1,0,O(x)+(\phi_0+O(x))t+O(t^2)\big)\\
            \Longrightarrow \text{Det}_t&=\frac{1}{x^2}\bigg(1+O(x^2)+(\phi_0O(x)+O(x^2))t+O(t^2)\bigg)\notag
        \end{aligned}
    \end{equation}
The additional \(\phi_0 O(x) t\) term destroys the "partially even" structure, causing the estimate for \(\partial_t III\) to fail in this case. In fact, the partially even structure is one of the key features underlying the theory of renormalized area. It plays a crucial role in ensuring the well-definedness of \(\operatorname{RenA}\), the existence of the universal lower bound in \textbf{Theorem~\ref{thm3}}, the estimate for term \(III\), and several other foundational results. In contrast, \cite{alexakis2010renormalized} employed a completely different approach, and therefore did not encounter this restriction. However, our results do not rely on the PE structure, making them applicable in a broader setting.
\end{remark}

From the theorems above, we know that a smooth critical point of \(\operatorname{RenA}\) must satisfy \(u_3 = 0\), therefore we have

\begin{corollary}
    Let $(M,g,V)$ satisfies \textbf{Condition C}. Suppose $h_3<0$ pointwise, then the infimum
    \begin{equation}
        I=\inf \{i(\Gamma):\Gamma \text{ closed smooth curves in } \Sigma\}
    \end{equation}
    can not be achieved by a smooth curve $\Gamma$.
\end{corollary}

\section{A Splitting Theorem}
In this section we are going to prove \textbf{Theorem \ref{thm3}}, i.e. show that given a smooth closed curve $\Gamma$, the infimum of renormalized area associated to $\Gamma$ can be achieved. Based on these minimizers, we then use \textbf{Theorem~\ref{thm1}} to prove a rigidity theorem when $(\Sigma,h)$ is a flat torus $\mathbb{T}^2$.

\subsection{Existence of Minimizer} ~

The proof will be divided into three parts: i) \(\operatorname{RenA}(F)\) is bounded below for \(F \in C(\Gamma)\); ii) Taking a minimizing sequence \(F_i\), we show that \(F_i\) subconverges to a minimal surface \(F_0\); iii) \(F_0\) is the surface we are seeking.

Before proceeding with the proof, let us recall the following theorem from geometric measure theory (e.g. \cite{Fdrer1960NormalAI}):

\begin{theorem}\label{GMT}
    Let $N^n$ a smooth manifold and $K\subset N$ compact. For each $c\in \mathbb{R}_+$, the set
    \begin{equation}
        \{\text{smooth manifold } F^{n-1}: \quad F\subset {K},\quad \mathcal{H}^{n-1}(F)+\mathcal{H}^{n-2}(\partial F)<c \}
    \end{equation}
    is precompact in weak topology in $\mathscr{I}^{n-1}$ where $\mathcal{H}^{p}$ is the $p$-dimensional Hausdorff measure and $\mathscr{I}^{n-1}$ is the set of integral $(n-1)$-currents.
\end{theorem}

\begin{proof}

\noindent \textit{i)Finiteness of $i(\Gamma_s)$}

 From \textbf{Lemma \ref{3_g_expansion}}, there exists a bdf so that
    \begin{equation}
        g=\frac{1}{x^2}(dx^2+ds^2+f^2d\theta^2+x^2\tilde{h}(r))
    \end{equation}
    in a cylinder region $D_{\epsilon_0}=\{x<\epsilon_0\}$, and $\tilde{h}(x)$ is a family of smooth metrics on $\Sigma$. The dependance on $x$ will be dropped for simplicity. Let $\Sigma_{\epsilon}=\{p:x(p)=\epsilon\}$. For $F\in C(\Gamma_s)$, again write
\begin{equation}
        RenA(F)=A(F\cap M_{\epsilon_0})-\frac{L(\Gamma_s)}{\epsilon_0}+\lim_{\epsilon\rightarrow 0}\bigg(A(F\cap(M_{\epsilon}\setminus M_{\epsilon_0}))-(\frac{1}{\epsilon}-\frac{1}{\epsilon_0})L(\Gamma_s))\bigg)\label{5_RenA}
        \end{equation}
It suffices to show that $\lim_{\epsilon\rightarrow 0}\bigg(A(F\cap(M_{\epsilon}\setminus M_{\epsilon_0}))-(\frac{1}{\epsilon}-\frac{1}{\epsilon_0}L(\Gamma_s))\bigg)$ is bounded below for arbitrary $F\in C(\Gamma)$.
     Suppose in $D_{\epsilon_0}$, $F$ can be expressed as a graph $u(x,\theta)$. $u=O(x^2)$ since $F\perp \Sigma$. Let $\partial_1=(1,0,\partial_xu)$ and $\partial_{2}=(0,1,\partial_\theta u$), then a similar calculation as in \textbf{Lemma \ref{lemma_det_expansion}} shows gives us $\text{Det}_F$. But unlike the computations in Section 4, we need to show that $C$ constants only depend on $g$ but not $F$, so we presents it here. In this the following $u_x=\partial_xu$ an $u_\theta=\partial_{\theta}u$.
    \begin{equation}
    \begin{aligned}
        g(\partial_1,\partial_1)&=\frac{1}{x^2}(1+u_x^2+x^2\tilde{h}_{ss}u_x^2)\\
        g(\partial_2,\partial_2)&=\frac{1}{x^2}(f^2+u_{\theta}^2+x^2\tilde{h}_{ss}u_{\theta}^2+2x^2\tilde{h}_{\theta s}u_{\theta})\\
        g(\partial_1,\partial_2)&=\frac{1}{x^2}(u_xu_{\theta}+x^2\tilde{h}_{ss}u_xu_{\theta}+x^2\tilde{h}_{\theta s}u_x)\notag
    \end{aligned}
    \end{equation}
    Direct computation shows that
    \begin{equation}
        \begin{aligned}
            |Det_F|^2&=\frac{1}{x^4}(g(\partial_1,\partial_1)g(\partial_2,\partial_2)-g(\partial_1,\partial_2)^2)\\
            &=\frac{1}{x^4}\bigg(f^2+f^2u_x^2+u_{\theta}^2+x^2(\tilde{h}_{ss}(f^2u_x^2+u_{\theta}^2)+2\tilde{h}_{\theta s}u_{\theta})-x^4\tilde{h}_{\theta s}u_x^2\bigg)\\
            &\geq \frac{1}{x^4}(f^2+\frac{1}{2}\min\{\inf f^2,1\}(u_x^2+u_{\theta}^2))\label{5_detF}
        \end{aligned}
    \end{equation}
if $x<\epsilon_0$ for a fixed small $\epsilon_0$. As a result
\begin{equation}
    \begin{aligned}
        &\quad \lim_{\epsilon\rightarrow 0}\bigg(A(F\cap(M_{\epsilon}\setminus M_{\epsilon_0}))-(\frac{1}{\epsilon}-\frac{1}{\epsilon_0}L(\Gamma_s)))\bigg)\\
        &=\lim_{\epsilon\rightarrow 0}\int_{\epsilon}^{\epsilon_0}\int_{S_{\theta}}\text{Det}_F-\frac{1}{x^2}d_\theta ds\geq 0\label{5-tail}
    \end{aligned}
\end{equation}
     Here we assumed $F$ is a graph in the region $D_{\epsilon_0}$. If it's not, we could extend the coordinates $(x,\theta,s)$ to the whole $\bar{M}$. Then consider the projection of $F$ into the plane $\{x,\theta,0\}$, denoted by $\Pi$, and apparently it's surjective. Apply Sard's theorem to $\Pi$, and the critical values has measure zero. For each regular value $p\in D_{\epsilon_0}$, we can find a neighborhood of $p$ in $F$ so that $F$ can be expressed locally as a graph, and the estimate $\text{Det}_F-\frac{1}{x^2}\geq 0$ still holds, so the same argument in (\ref{5-tail}) works if we replace the first equality into inequality because $\Pi$ might be a multiple-to-one map somewhere.

\noindent \textit{ii):Compactness}

Since $i(\Gamma_s)$ is finite, we can find a minimizing sequence $\{F_i\}$ for $i(\Gamma_s)$. Let $a_n=\frac{1}{e^n}$. We are going to use \textbf{Theorem \ref{GMT}} to show that for each fixed $n\in\mathbb{Z}_+$, we can find $\epsilon_n\in (a_{n+1},a_n)$ so that $\{F_i\}$ subconverge to an area-minimizing smooth manifold in $M_{\epsilon_n}=\{x\geq \epsilon_n\}$. To show this, we need to find $\epsilon_n\in(a_{n+1},a_n)$ and $c_n$ such that
\begin{equation}
    \mathcal{H}^{n-1}(F_i\cap M_{\epsilon_n})+\mathcal{H}^{n-2}(\partial (F_i\cap M_{\epsilon_0}))<c_n
\end{equation}
for all $i$.

The idea is, even though for $L(F_i\cap M_{\epsilon})$ can be large for small $\epsilon$, but it can not be large in an interval, so we can find $\epsilon_n\in (a_{n+1},a_n)$ so that $L(F_i\cap M_{\epsilon_n})$ is bounded.

Fix $n$. Let $E_n$ be the annular region $\{a_{n+1}<x<a_n\}$. Since $RenA(F_i)$ is bounded from above, from (\ref{5_RenA})(\ref{5-tail}) we know that $A(F\cap E_{n})$ is bounded above. Assume $F_i$ is a graph $u(x,\theta)$. Then following the second line of (\ref{5_detF})
\begin{equation}
    \begin{aligned}
        A(F_i\cap E_n)&\geq \int_{a_n+1}^{a_n}\int_{S_{\theta}}\frac{1}{x^2}\bigg(1+\frac{1}{2}\min\{\inf f^2,1\}(u_x^2+u_{\theta}^2)-Cx^2\bigg)d\theta dx\\
        &\geq \frac{1}{2}\min\{\inf f^2,1\}\int_{a_n+1}^{a_n}\int_{S_{\theta}}\frac{1}{x^2}u_{\theta}^2d\theta dx+C(a_n,g)\label{5-1-middle1}\\ 
    \end{aligned}
\end{equation}
Here \textbf{$C(a_n,g)$ is some bounded constant that depends on $a_n$ and $g$ but not $i$}. Consider the closed curve $\sigma_{i,x_0}=F_i\cap \Sigma_{x_0}=\{(x_0,\theta,u(x_0,\theta))\}$ for $x\in (a_{n+1},a_n)$. Its tangent vectors are $\partial_{\theta}=(0,1,u_{\theta})$ and thus its length is controlled by
\begin{equation}
\begin{aligned}
    L(\sigma_{i,x},g)&=\int_{\sigma_{i,x}}\frac{1}{x}\bigg(1+u_{\theta}^2+x^2(\tilde{h}_{ss}u_{\theta}^2+2\tilde{h}_{\theta s}u_{\theta})\bigg)^{\frac{1}{2}}d\theta\\
    &\leq \int_{\sigma_{i,x}} \frac{1}{x}(1+u_{\theta}^2+C(g)x^2)\\
    &\leq \int_{\sigma_{i,x}}\frac{1}{x}u_{\theta}^2d\theta+C(a_n,g)\notag
\end{aligned}
\end{equation}
So (\ref{5-1-middle1}) becomes
\begin{equation}
    A(F_i\cap E_n)\geq \frac{1}{2}\min\{\inf f^2,1\}\int_{a_{n+1}}^{a_n}\frac{1}{x}L(\sigma_{i,x},g)dx+C(a_n,g)\label{5-1-middle2}
\end{equation}

If \(F_i\) is not graphical over the \((s, \theta)\)-coordinates, we can first perturb \(F_i\) slightly so that each slice \(\sigma_{i,x}\) consists of a finite number of graphical curves and finitely many isolated points for each \(x \in (a_{n+1}, a_n)\), without significantly altering the area. (Note that in some cases, \(F_i\) may be vertical, so that \(\sigma_{i,x}\) is a 2-dimensional surface rather than a curve.) After this perturbation, we can apply the argument above locally near each graphical curve to obtain $(\ref{5-1-middle1})$.

Since \(F_i\) is a minimizing sequence for $RenA$, it follows from equations~\eqref{5_RenA} and \eqref{5-tail} that the area \(A(F_i \cap M_{a_{n+1}})\) is uniformly bounded from above in \(i\) (though it may depend on \(n\)). Consequently, the area \(A(F_i \cap E_{a_{n+1}})\) is also uniformly bounded.
 Then from (\ref{5-1-middle2}),
\begin{equation}
    \int_{a_{n+1}}^{a_n}\frac{1}{x}L(\sigma_{i,x},g)dx<c(n,g)\label{5-1-middle3}
\end{equation}
where $c(n,g)$ is a constant depending on $n,g$ and fixed in the following claim.

\noindent\textit{Claim: Let interval $I_{i,n}\coloneqq \{x\in (a_{n+1},a_n): L(\sigma_{i,x},g)<2c(n,g)\}$, then $\mathcal{L}({I}_{i,n})>c'>0$ where $c'$ is independant of $i$.}

\noindent\textit{Proof of claim}\\

If not, let $I^c_{i,n}=(a_{n+1},a_n)\setminus I_{i,n}$, then following (\ref{5-1-middle3})
\begin{equation}
\begin{aligned}
    c(n,g)>\int_{a_{n+1}}^{a_n}\frac{1}{x}L(\sigma_{i,x},g)dx&\geq \int_{I_{i,n}^c}\frac{1}{x}L(\sigma_{i,x},g)dx\\
   & \geq 2c(n,g)\int_{a_{n+1}}^{a_{n+1}+\mathcal{L}(I_{i,n}^c)}\frac{1}{x}dx\\
   &\geq 2c(n,g)\big(\log (a_{n+1}+\mathcal{L}(I_{i,n}^c))-\log (a_{n+1})\big)\notag
\end{aligned}
\end{equation}
If \(\mathcal{L}(I_{i,n}^c) = a_n - a_{n+1}\), recall \(a_n = e^{-n}\) and we are led to a contradiction. Therefore, \(\mathcal{L}(I_{i,n}^c)\) must be bounded above by a constant strictly less than \(a_n - a_{n+1}\). It follows that \(\mathcal{L}(I_{i,n}) > c' > 0\) for some constant \(c'\) independent of \(i\).
\qed

Now since $\mathcal{L}(\cup_{i=k}^{\infty}I_{i,n})\leq a_{n}-a_{n+1}<\infty$ for each $k$, we have by the continuity of Lebesgue measure that
\begin{equation}
   \mathcal{L}\big( \lim_{k\rightarrow \infty} (\cup_{i=k}^{\infty}I_{i,n})\big)=\lim_{k\rightarrow \infty} \mathcal{L}\big(\cup_{i=k}^{\infty}I_{i,n}\big)>c'
\end{equation}
Therefore, the union \(\bigcup_{i = k}^{\infty} I_{i,n}\) is non-empty, and we choose \(\epsilon_n\) to be an element in this set. By definition, there exists a subsequence (still denoted by \(\{F_i\}\)) such that $A(F_i \cap M_{\epsilon_n}) + L(F_i \cap \Sigma_{\epsilon_n})$ is uniformly bounded in \(i\). Then, by \textbf{Theorem~\ref{GMT}}, the sequence \(\{F_i\}\) subconverges in the region $M_{\epsilon_n}$. Then let $n\rightarrow \infty$ and by a diagonalizing argument, we get a subsequence that converges to $F_0$. to Since \(\{F_i\}\) is a minimizing sequence for the renormalized area, it follows that \(F_0\) is an area-minimizing surface.

\noindent\textit{iii) Regularity and Infimum}

The interior regularity is guaranteed by Allard's regularity theorem \cite{Allard1972OnTF}, and boundary regularity is from \cite{lin1989dirichlet,hardt1987regularity,Tonegawa1996ExistenceAR}, and see \cite{alexakis2010renormalized} for a summary.

It remains to check that $F_0$ is a minimizer for $i(\Gamma_s)$. Since $F_0$ is regular to the boundary, thus can be expressed as as a graph $u_0$. With this $u_0$ fixed, we could run the calculation (\ref{5_detF})(\ref{5-tail}) to get that the following
\begin{equation}
    \lim_{\epsilon\rightarrow 0}\bigg(A(F_0\cap(M_{\epsilon}\setminus M_{\epsilon_0}))-(\frac{1}{\epsilon}-\frac{1}{\epsilon_0})L(\Gamma_s))<C\epsilon_0\label{5-1-middle4}
\end{equation}
where $C$ depends on $g$ and $F_0$, but not $\epsilon_1$.
To estimate \( RenA(F_0) -RenA(F_i) \), recall equation~\eqref{5_RenA}. For each $\delta$, by~\eqref{5-1-middle4} and~\eqref{5-tail} we can find $\epsilon_0$ small so that 
\begin{equation}
    \begin{aligned}
        &\lim_{\epsilon \to 0} \left( A(F_0 \cap (M_\epsilon \setminus M_{\epsilon_0})) - \left( \frac{1}{\epsilon} - \frac{1}{\epsilon_0} L(\Gamma_s) \right) \right)\\
        -&\lim_{\epsilon \to 0} \left( A(F_i \cap (M_\epsilon \setminus M_{\epsilon_0})) - \left( \frac{1}{\epsilon} - \frac{1}{\epsilon_0} L(\Gamma_s) \right) \right)\\
        <&\frac{\delta}{2}\label{5-1-middle5}
    \end{aligned}
\end{equation}
for all $i$. With \(\epsilon_0\) fixed, the sequence \(F_i\) converges to \(F_0\) in \(M_{\epsilon_0}\), i.e., $
A(F_0 \cap M_{\epsilon_0}) - A(F_i \cap M_{\epsilon_0}) < \frac{\delta}{2}$ for sufficiently large \(i\). Combine \eqref{5-1-middle5}, we have \( RenA(F_0) < RenA(F_i) + \delta \) for these values of \(i\). Since \(F_i\) is a minimizing sequence for \(i(\Gamma_s)\), it follows that \(F_0\) is a minimizer for \(i(\Gamma_s)\).

\end{proof}

\subsection{A Rigidity Result }

\begin{proof}
    Let $\phi(s)=i(\Gamma_s)$ for $s\in S_s$. For fixed $s_0$, by \textbf{Theorem \ref{thm3}}, there exists a $F_{s_0}$ minimizes $i(\Gamma_{s_0})$. Run the flow $\partial\Phi(\cdot,t)=V\nu$ from $F_{s_0}$ in $t\in [0,T]$ for a small number $T$, and let $F_t=\Phi(\cdot,t)$, $\psi(t)=RenA(F_t)$. Since $\Phi(\cdot,t)$ restricts to a normal flow on $\Sigma$ and $(\Sigma,h)$ is a flat torus, then $F_t\cap \Sigma=\Gamma_{s_0+t}$, i.e. $F_t\in C(\Gamma_{s_0+t})$. Therefore $\psi$ is a support function for $\phi(s)$ at $s_0$ from above. Since $\kappa=0$, from \textbf{Theorem \ref{thm2}}, we have
    \begin{equation}
        \frac{\partial^2}{\partial t^2}\phi\bigg|_{s_0}\leq \int_{\Gamma_{s_0}}\frac{3}{4}(tr_{h}h_3+h_3(\partial s,\partial s))d\theta_h
    \end{equation}
where the inequality is understood in the support sense. Since $\phi(s)$ is a periodic function, case (i) can not happen.

Now, suppose case (ii) occurs. In this case, \(\phi\) is a concave function, and thus \(\phi\) must be a constant function by its periodicity. We will show that 
\begin{equation}
    \frac{\partial}{\partial t}\psi(t_0)\leq 0,\quad 0<t_0<T\label{5-2-first}
\end{equation}
Since \(\phi\) is constant and \(\phi \leq \psi\), it follows that \(\psi\) must also be constant. This constancy of \(\psi\) implies a local splitting structure. The proof will be like (ii) of \textbf{Theorem \ref{thm2}}, but here $t_0$ is a fixed small number, while in  \textbf{Theorem \ref{thm2}} we are taking derivative so $t\rightarrow 0$. In this proof, \textbf{$C$ denotes some bounded function that doesn't depend on $t$}.

WLOG, calculate at $s_0=0$. Let $F_t$ be as above and near conformal boundary $F_0$ is a graph $u(x,\theta)=u_3x^3+O(x^4)$. Since $\psi(t)$ supports $\phi(t)$ from above at $0$ and $\phi$ is constant, 
\begin{equation}
    0=\frac{\partial}{\partial t}\phi(0)=\frac{\partial}{\partial t}\psi(0)=-\int_{\Gamma_0}3u_3d\theta\label{5-2-u3}
\end{equation}

Let $\gamma$ be a geodesic from $p\in \Sigma_{\epsilon}$ normal to $F_0$. We want to estimate $\dot{\gamma}(t_0)$ for a small fixed $t_0$. By Jacobi equation $\dot{\gamma}^{k}(t)=-\gamma^i(t)\gamma^j(t)\Gamma_{ij}^k(t)$ where $\Gamma_{ij}^k$ is Christoffel symbol for $\bar{g}$, and from (\ref{3_g_expansion})  and the assumption $h=d\theta^2+ds^2$ we have

\begin{equation}
    \bar{g}=dx^2+d\theta^2+ds^2+(h_3+(tr_{h}h_3)h)x^3+O(x^4)\label{5-2-g-expansion}
\end{equation}

Therefore $\Gamma_{ij}^k(\gamma(t))=C\epsilon^2$ for all $i,j,k$. $\dot{\gamma}(0)$ can be calculated to be
\begin{equation}
    \dot{\gamma}(0)=(-3u_3\epsilon^2+C\epsilon^3,C\epsilon^3,1+C\epsilon^3) \label{5-2-norm0}    
\end{equation}
 we get that $\dot{\gamma}(t)=(C\epsilon^2,C\epsilon^2,1+C\epsilon^2)$. It follows that
\begin{equation}
\begin{aligned}
    \dot{\gamma}^x(t_0)&=\dot{\gamma}^x(0)-\int_0^{t_0}\gamma^i(t)\gamma^j(t)\Gamma_{ij}^k(t)dt\\
    &=-3u_3\epsilon^2+\frac{3}{4}\int_0^{t_0}(h_{3,ss}+tr_{h_0}h_3)(\gamma(t))x^2dt+O(\epsilon^3)\\
    &=-3u_3\epsilon^2+\frac{3}{4}\epsilon^2\int_0^{t_0}(h_{3,ss}+tr_{h_0}h_3)(\gamma(t))dt+O(\epsilon^3)\\
    \dot{\gamma}^{\theta}(t_0)&=C\epsilon^3\\
    \dot{\gamma}^{s}(t_0)&=1+C\epsilon^3
\end{aligned}    
\end{equation}
The second and third equality is calculated in a similar way using $\Gamma^{\theta}_{ss}=\Gamma^{s}_{ss}=C\epsilon^3$.

Since we are running normal flow, $\dot{\gamma}$ is parallel to $\bar{\nu}$. After normalization and an argument similar to \textbf{Lemma \ref{lemma_w}} which deals with shifting in $\theta$ and $x$, we see that on $F_{t_0}$
\begin{equation}
    \bar{\nu}=(-3u_3\epsilon^2+\frac{3}{4}\epsilon^2\int_0^{t_0}(h_{3,ss}+tr_{h_0}h_3)(\gamma(t))dt,C\epsilon^3,1+C\epsilon^3)\label{5-2-norm}
\end{equation}
Using equation~\eqref{4_initial}, we again split \(\operatorname{RenA}(F_t)\) into three parts: \(I\), \(II\), and \(III\). Combining this with equation~\eqref{5-2-norm} and the fact that \(L(\Gamma_s)\) remains constant, and applying the same strategy as in \textit{ii)} of the proof of \textbf{Theorem~\ref{thm2}}, we obtain

\begin{equation}
\begin{aligned}
    \frac{\partial}{\partial t}II(t_0)&=\int_{S_{\theta}}-3u_3d\theta+\frac{3}{4}\int_0^{t_0}\int_{S_{\theta}}(h_{3,ss}+tr_{h_0}h_3)d\theta ds+O(\epsilon)\\
    &\leq \int_{S_{\theta}}-3u_3d\theta+O(\epsilon)=O(\epsilon)\label{5-2-II}
\end{aligned}
\end{equation}
where the inequality in the second line follows from the assumption of the theorem and the last equality follows from (\ref{5-2-u3}). $\partial_t I$ can be computed as
\begin{equation}
    \begin{aligned}
        \frac{\partial}{\partial t}I(t_0)=\int_{F_{t_0}\cap M_\epsilon}Hdv_g\leq 0\label{5-2-I}
    \end{aligned}
\end{equation}
where the second inequality follows from equation~\eqref{4-dHdt} by applying Grönwall's inequality and using the fact that \(H = 0\) on \(F_0\). From (\ref{5-2-norm}), there is no $\epsilon$ term in $x$ direction, and $\partial_t L(\Gamma_t,h)=0$, so the same argument in (\ref{4_dIIIdt}) shows that $\partial_tIII=C\epsilon$. Together with (\ref{5-2-II})(\ref{5-2-I}), let $\epsilon\rightarrow 0$, we have
\begin{equation}
    \frac{\partial}{\partial t}\psi(t)\leq \int_{F_{t_0}}Hdv_g\leq 0
\end{equation}
Therefore, \(\psi\) must be constant, since it supports the constant function \(\phi\) from above. It follows that \(H = 0\) for each surface \(F_t\). Then, from equation~\eqref{4-dHdt}, we conclude that \(b = 0\) on each \(F_t\). As a result, the metric \(g\) locally splits as
\[
g = V^2 dt^2 + \tilde{g},\label{5-2-splitting}
\]
with \(t = s\) on \(\Sigma\). The coefficient \(V^2\) arises from the fact that \(\partial_t = V \nu\), so that \(|\partial_t|_g = V\). And by the splitting, the flow $F_t$ can be extended to $\mathbb{R}$.

Next, we show that \(V\) is independent of \(t\). Most of the computation can be found in \cite{GSW, wang2024riccati}; we include the entire argument here for completeness. In the following, all covariant derivatives and curvature terms are taken with respect to the metric \(g\) unless specified otherwise. Let \(p,q\) denote indices corresponding to a local normal coordinate along $F_t$.
 Direct computation shows that 
\begin{equation}\label{5-2-connection}
        \nabla_{\partial t}\partial t=-VV_p\partial_p+\frac{V_t}{V}\partial t;\quad \nabla_{\partial_p}\partial_t=\nabla_{\partial_t}\partial_p=\frac{V_p}{V}\partial t;\quad \nabla_{\partial_p}\partial_q=0.
    \end{equation}
Since $F_t$ are totally geodesic, by Codazzi equation we know that $Ric_{g}(\partial_p,\partial_t)=0$. It follows from (\ref{static_Einstein_equ}) that $\nabla^2V(\partial _t,\partial_p)=0$, and combine with (\ref{5-2-connection}) we get
         \begin{equation}
     \begin{aligned}
         0=\nabla^2V(\partial_t,\partial_p)&=\partial_p\partial_tV-(\nabla_{\partial_t}\partial_p)V\\
         &=\partial_p\partial_tV-\frac{\partial_pV\partial_tV}{V}=V\partial_p\partial_t \log(V)\notag
     \end{aligned}
     \end{equation}
     which implies the splitting $V(t,y)=\alpha(t)\beta(y)$ for $y$ coordinate for $F_t$. Now examine (\ref{static_Einstein_equ}) along $F_t$ direction. From (\ref{5-2-connection}) we can compute
    \begin{equation}
        \begin{aligned}
             R(\partial_p,\partial t,\partial t,\partial_p)=\langle\nabla_{\partial t}\nabla_{\partial p}\partial_p,\partial t\rangle-\langle\nabla_{\partial p}\nabla_{\partial r}\partial_p,\partial t\rangle=-\alpha^2\beta\beta_{ii}
             \label{3_computation_curvature_1}
        \end{aligned}
    \end{equation}
    Using Gauss equation we have
    \begin{equation}
       Ric(\partial_p,\partial_p)=Ric_{\tilde{g}}(\partial_p,\partial_p)+R(\partial_p,\frac{\partial t}{V},\frac{\partial t}{V},\partial_p)=Ric_{\tilde{g}}(\partial_p,\partial_p)-\frac{\beta_{ii}}{\beta}\label{4_Ric_computation}
    \end{equation}
    And $\nabla^2V$ can be computed as
    \begin{equation}
    \begin{aligned}\label{4_hessian}
        \nabla^2V(\partial_p,\partial_p)&=\alpha\beta_{pp}\\
        \nabla^2V(\frac{\partial_t}{V},\frac{\partial_t}{V})
        &=\frac{\alpha_{rr}}{\alpha^2\beta}+\frac{\alpha|\nabla_{h}\beta|^2}{\beta}-\frac{\alpha_r^2}{\alpha^3\beta}\\        
        \Rightarrow \Delta V&=\alpha(\Delta_{h}\beta+\frac{|\nabla_{\tilde{g}}\beta|^2}{\beta})+\frac{1}{\beta}(\frac{\alpha_{rr}}{\alpha^2}-\frac{\alpha_r^2}{\alpha^3})
        \end{aligned}
    \end{equation}
    Put these together, we get
    \begin{equation}
    \begin{aligned}
        0&=\sum_{p}\big(VRic(\partial_p,\partial_p)+\Delta V-\nabla^2V(\partial_p,\partial_p)\big)\\
        &=\alpha\big[\beta R_{\tilde{g}}+2\frac{|\nabla_{\tilde{g}}\beta|^2}{\beta})\big]+\frac{2}{\beta}\big[\frac{\alpha_{rr}}{\alpha^2}-\frac{\alpha_r^2}{\alpha^3}\big]\\
       \Rightarrow & \beta\big[\beta R_{\tilde{g}}+2\frac{|\nabla_{\tilde{g}}\beta|^2}{\beta})\big]=-\frac{2}{\alpha}\big[\frac{\alpha_{rr}}{\alpha^2}-\frac{\alpha_r^2}{\alpha^3}\big]
    \end{aligned}
    \end{equation}
    This equation is seperable, which implies that $\frac{1}{\alpha}\big[\frac{\alpha_{rr}}{\alpha^2}-\frac{\alpha_r^2}{\alpha^3}\big]=C$. If we make the coordinate transformation $dt=\alpha dr$, this implies that
    \begin{equation}
        \alpha \frac{\partial^2 \alpha}{\partial t^2}=C \label{3_middle_4}
    \end{equation}
    for some constant $c$. The only solution bounded below is the constant function, so $V$ does not depend on $t$.
Finally, the rigidity for Horowitz-Myers solition comes from \textbf{Theorem IV} in \cite{GSW}.
\end{proof}

In \cite{GSW, wang2024riccati}, Busemann function method is applied to prove a similar result:
\begin{theorem}[G.Galloway, S.Surya, E.Woolgar]\label{thm_splitting}
    Consider $(M,g,V)$ satisfying \textbf{Condition C} We also assume that: the second fundamental form $b$ for level sets $\{\frac{1}{V}=\epsilon\}$ is positive semi-definite for small $\epsilon$. Then the Riemannian universal cover $(\tilde{M}^*,\tilde{g}^*)$ of $(\tilde{M},\tilde{g}=\frac{1}{V^2}g)$ splits isometrically as
    \begin{equation}
        \tilde{M}^*=\mathbb{R}^k\times\Sigma,\quad \tilde{g}^*=g_{\mathbb{E}}+ \tilde{h}
    \end{equation}
    where $(\mathbb{R}^k,g_{\mathbb{E}})$ is standard $k$-dimensional Euclidean space with $0\leq k\leq n$, and $(\Sigma,\tilde{h})$ is a compact Riemannian manifold with non-empty boundary. Furthermore, both $h$ and $V^*$ are independent of $r$. As a result, the Riemannian universal cover $(M^*,g^*)$ of $(M,g)$ splits isometrically as
    \begin{equation}
        M^*=\mathbb{R}^k\times \Sigma,\quad g^*=(V^{*2}g_{\mathbb{E}})+ h
    \end{equation}
    where $V^*=V\circ \pi$($\pi=$ covering map) is constant along $\mathbb{R}^k$.
\end{theorem}

In this theorem, rigidity was obtained by assuming that the second fundamental form of the level sets \(\Sigma_{\epsilon} = \left\{ \frac{1}{V} = \epsilon \right\}\) is semi-convex with respect to the metric \(\bar{g} = \frac{1}{V^2} g\) for small \(\epsilon\). If \(h'_3\coloneqq h_3+(tr_hh_3)h < 0\) on \(\Sigma\), then this convexity condition is clearly satisfied. However, if we only assume \(h'_3 \leq 0\), difficulties arise. 

For example, suppose \(h'_3\) is taken to be that of the Horowitz–Myers soliton. In this case, all higher-order coefficients \(h'_i\) for \(i \geq 4\) in the expansion
\[
\frac{1}{V^2} g = h + h'_3 x^3 + \sum_{i \geq 4} h'_i x^i
\]
are uniquely determined by \(h\) and \(h'_3\) (we use ' to distinguish it from the expansion for $x^2g$) , and therefore coincide with those of the Horowitz–Myers soliton. In particular, we have
\[
h'_i(\partial_s, \partial_s) = 0, \quad \forall i \geq 4.
\]
In such a scenario, the convexity of the level sets \(\Sigma_\epsilon\) cannot be determined solely from \(h_3\), suggesting that additional assumptions may be necessary.

By contrast, our result is entirely based on \(h'_3\) and avoids the need for convexity assumptions on the level sets. Besides, our approach only require an assumption on the integral of \(h'_3\) (rather than a pointwise condition). However, our assumption is more intricate and dimensionally restricted.

Moreover, our method also applies to warped product structures. By the Gauss--Bonnet theorem, we have \(\int_{\Sigma} R_h = 0\). Therefore, from the expansion in equation~\eqref{3_g_expansion}, the assumption that the level sets \(\Sigma_\epsilon\) are semi-convex is only possible if \((\Sigma, h)\) is flat. In contrast, our approach remains valid for warped product metrics. In particular, the warped product structure is used to ensure that \(F_t \in C(\Gamma_{s_0 + t})\), so that \(\psi\) serves as a support function for \(\phi\); see the first paragraph of the proof for details.\\

As discussed above, our argument also requires a convexity assumption, but in an integral sense. This idea is illustrated in Figure~\ref{renA4}. For simplicity, we assume that \(h_3 < 0\) pointwise.

Following the proof of \textbf{Theorem~\ref{thm2}}, without loss of generality, assume that the infimum of \(\phi(s) = i(\Gamma_s)\) is achieved at \(s = 0\). Let \(F_0\) be a minimizer for \(i(\Gamma_0)\), and suppose \(F_0\) is represented as a graph \(u = u_3 x^3 + O(x^4)\). Since \(s = 0\) is a critical point for \(\phi\), we have $\int_{\Gamma_0} u_3 \, d\theta = 0.$

In Figure~\ref{renA4}, \(\gamma\) denotes a geodesic emanating from \(\Sigma_\epsilon \cap F_0\) and normal to \(F_0\). Because \(\Sigma_\epsilon\) is convex, \(\gamma\) deviates to the right of \(\Sigma_\epsilon\), so we observe from the figure that
\[
F_{t,\epsilon} \subsetneq \Phi(F_{0,\epsilon}, t).
\]
This illustrates how the sign of \(h_3\) influences the second variation of the renormalized area.

Recall that, starting from a minimal surface, the area is non-increasing along the flow \(\partial_t \Phi = V \nu\). Therefore, we have
\begin{equation}
    A(F_{t,\epsilon}) < A(\Phi(F_{0,\epsilon}, t)) \leq A(F_{0,\epsilon}).
\end{equation}
Since $L(\Gamma_t)=L(\Gamma_0)$, taking $\epsilon\rightarrow 0$, we arrive at $RenA(F_t)< RenA(F_0)$, establishing the instability.\\

One might ask whether these results can be generalized to higher dimensions. In even dimensions, the renormalized area is not well-defined due to the appearance of anomalies, which arise from different choices of boundary defining functions (bdf's). If we assume \textbf{Condition C}, then \(\frac{1}{V}\) provides a canonical choice of bdf, making \(\operatorname{RenA}\) well-defined. However, logarithmic terms may still appear in the expansion, which may obstruct the generalization.

In odd dimensions, a generalization is more promising, although the calculations become significantly more involved. If we assume that the boundary \((\Sigma, h)\) has vanishing Ricci curvature, then for an asymptotically Poincar\'e--Einstein (APE) manifold \((M^n, g)\), there exists a bdf \(x\) such that
\[
g = \frac{1}{x^2} \left( dx^2 + h + x^{n-1} h_{n-1} + O(x^n) \right),
\]
and in this setting, \textbf{Theorem~\ref{thm1}} should carry over to higher dimensions. Consequently, \textbf{Theorem~\ref{thm2}} is also expected to admit a natural generalization.

\textbf{Acknowledgements}

The author is thankful to Prof Biao Ma, Jared Marx-Kuo, Jeffrey S. Case, Rafe Mazzeo and for their kind help. 

\end{document}